\documentclass[11pt]{amsart}
\usepackage{amsmath,amssymb,amsfonts,url,mathptmx}
\numberwithin{equation}{section}

\newtheorem{thm}{Theorem}[section]
\newtheorem{lem}{Lemma}[section]
\newtheorem{rem}{Remark}[section]
\newtheorem{prop}{Proposition}[section]

\begin{document}
\title[Singular Liouville Equation]{On Non-simple blowup solutions of Liouville equation}

\author{Teresa D'Aprile}
\address{Dipartimento di Matematica, Universit\`a di Roma “Tor Vergata”, via della Ricerca
Scientifica 1, 00133 Roma . Italy.} \email{daprile@mat.uniroma2.it}

\author{Juncheng Wei}
\address{Department of Mathematics \\ University of British Columbia\\ Vancouver, BC V6T1Z2, Canada} \email{jcwei@math.ubc.ca }

\author{Lei Zhang} \footnote{The research of J. Wei is partially supported by NSERC of Canada. Lei Zhang is partially supported by a Simons Foundation Collaboration Grant}
\address{Department of Mathematics\\
        University of Florida\\
        1400 Stadium Rd\\
        Gainesville FL 32611}
\email{leizhang@ufl.edu}

\date{\today}

\begin{abstract} For Liouville equation with quantized singular sources, the non-simple blowup phenomenon has been a major difficulty for years. It was conjectured by the first two authors that the non-simple blowup phenomenon does not occur if the equation is defined on the unit ball with Dirichlet boundary condition. In this article we not only completely settle this conjecture in its entirety, but also extend our result to cover any bounded domain. Since the main theorem in this article rules out the non-simple phenomenon in commonly observed applications, it may pave the way for advances in degree counting programs, uniqueness of blowup solutions and construction of solutions, etc.
\end{abstract}

\maketitle

\section{Introduction}

In this article we study the following Liouville equation with finite singular sources: let $\Omega$ be an open, bounded and connected subset of $\mathbb R^2$ with smooth boundary $\partial \Omega$ and let $u$ be defined as
\begin{equation}\label{main-e-1}
\Delta u+e^u=\sum_{t=1}^M4\pi \gamma_t \delta_{p_t} \quad \mbox{in}\quad \Omega
\end{equation}
where $p_1,...,p_M$ are points in $\Omega$, $\gamma_1,...,\gamma_M$ are constants greater than $-1$. If a $\gamma_t$ is a positive integer ($\gamma_t\in \mathbb N$, $\mathbb N$ is the set of natural numbers), we say $p_t$ is a \emph{ quantized singular source }. Studying the asymptotic behavior of blowup solutions 
is particularly important for a number of applications.    Let $\mathfrak{u}_k$ be a sequence of solutions of (\ref{main-e-1}); we say $\{\mathfrak{u}_k\}$ is a sequence of blowup solutions if, for a point $p\in \Omega$, there exist $x_k\to p$ such that $\mathfrak{u}_k(x_k)-2\gamma_p\log |x_k-p|\to \infty$. Here we set $\gamma_p=0$ if $p$ is not a singular source. For a sequence of blowup solutions $\mathfrak{u}_k$, it is standard to assume a uniform bound on the total integration and boundary oscillation: there exists $C>0$ independent of $k$ such that 

\begin{equation}\label{standard-a}
\int_{\Omega}e^{\mathfrak{u}_k}\le C
\end{equation}
and 
\begin{equation}\label{standard-b}
|\mathfrak{u}_k(x)-\mathfrak{u}_k(y)|\le C, \quad \forall x,y\in \partial \Omega. 
\end{equation}

If a blowup point $p$ is either a regular point or a ``non-quantized" singular source, the asymptotic behavior of $\mathfrak{u}_k$ around $p$  is well understood (see \cite{bclt,bart4,chen-lin-cpam-02,chen-lin-cpam-15,gluck,li-cmp,li-shafrir,malchiodi,zhangcmp,zhangccm}). As a matter of fact, $\mathfrak{u}_k$ satisfies the spherical Harnack inequality around $p$, which implies that, after scaling, the sequence $\mathfrak{u}_k$  behaves as a single bubble around  the maximum point. However, if $p$ happens to be a quantized singular source, the so-called ``non-simple" blowup phenomenon does happen (see \cite{kuo-lin-jdg, wei-zhang-adv,wei-zhang-plms,wei-zhang-jems}), which is equivalent to stating that $\mathfrak{u}_k$ violates the spherical Harnack inequality around $p$. The study of non-simple blowup solutions has been a major challenge for Liouville equations and its research laid dormant for years. Recently significant progress has been made by Kuo-Lin, Bartolucci-Tarantello and other authors \cite{bart3,bart-taran-jde-2,daprile-wei-jfa,kuo-lin-jdg,wei-zhang-adv,wei-zhang-plms,wei-zhang-jems}. In particular  it is established in \cite{bart-taran-jde-2} and \cite{kuo-lin-jdg} that  there are $\gamma_p + 1$ local maximum points and they are evenly distributed on $\mathbb{S}^1$ after scaling according to their magnitude.

 The case $\gamma_p\in\mathbb{N}$ is more difficult to treat, and at the same time the most relevant to physical applications. Indeed, in vortex theory the number
$\gamma_p$
represents vortex multiplicity, so that  in that context the most interesting case is
precisely
when it is a positive integer. The difference between the case $\gamma_p\in \mathbb{N}$ and $\gamma_p\not\in \mathbb{N}$ is also analytically essential. Indeed,  as usual in
problems involving  concentration phenomena like \eqref{main-e-1},  after
suitable rescaling of the blowing-up around a concentration point one sees a limiting equation which, in this case, takes the form of the planar singular Liouville equation:
$$\Delta U+e^U=4\pi\gamma \delta_p\;\;\hbox{ in } \mathbb{R}^2, \quad \int_{\mathbb R^2} e^U dx <\infty;$$ 
 only if $\gamma_p\in\mathbb{N}$ the above limiting equation   admits non-radial solutions around $p$
 since the family of all solutions extends to one carrying an extra parameter  (see \cite{prajapat}). This suggests that if $\gamma_p\in\mathbb{N}$  and the blow-up point happens to be the singular source, then  solutions of \eqref{main-e-1}  may exhibit non-simple blow-up phenomenon. 

So, from analytical viewpoints the study of non-simple blowup solutions is far more challenging than simple blowup solutions, but the impact of this study may be even more significant because non-simple blowup solutions represent certain situations in the blowup analysis of systems. If local maximums of blowup solutions in a system tend to one point, the profile of solutions can be described by a Liouville equation with quantized singular source. It is desirable to know exactly when non-simple blowup phenomenon happens. In \cite{daprile-wei-jfa} the first two authors studied the following equation:
      \begin{align} \label{dir-0-cond}
\Delta u+\lambda e^u&=\sum_{t=1}^M 4\pi \gamma_t\delta_{p_t} \quad \mbox{in}\quad \Omega\subset \mathbb R^2, \\
u&=0 \quad \mbox{on}\quad \partial \Omega, \nonumber
\end{align}
where $\Omega$ is an open and bounded subset of $\mathbb R^2$, $p_1,...,p_M\in \Omega$, $\partial \Omega$ is smooth, $\lambda>0$ and $\gamma_t\in \mathbb N$. We note equation \eqref{dir-0-cond} transforms into a Liouville equation with no boundary condition of the type \eqref{main-e-1} by using the presence of the free parameter $\lambda>0$ under the transformation $\mathfrak{u}= u(x)+\log \lambda$; then, we say that $u_k$ is a sequence of blow-up solutions of  (\ref{dir-0-cond}) with parameter $\lambda_k\to 0^+$ if, for some point $p\in\Omega$, there exists $x_k\to p$ such that  
$u_k(x)+\log \lambda_k-2\gamma_p\log|x_k-p|\to +\infty$ 
and there is a uniform  bound on its mass:
\begin{equation}\label{unif-bound-3}
\lambda_k\int_{\Omega}e^{u_k}<C. 
\end{equation}
In \cite{daprile-wei-jfa}
it was 
conjectured that if $\Omega=B_1$, $M=1$ and $p_1=0$, then there is no non-simple blowup sequence in $B_1$.  
This conjecture of D'Aprile and Wei was considered audacious since there was no restriction of number of bubbles inside $\Omega$ and only the Dirichlet boundary condition is placed. The purpose of this conjecture is to claim that the boundary data seem to have a great influence on the profile of blowup solutions inside. In comparison, the series of works of Wei-Zhang \cite{wei-zhang-adv,wei-zhang-plms,wei-zhang-jems} focuses on the vanishing rate of coefficient functions. In this respect, the conjecture of D'Aprile and Wei seems more natural and useful for application. In the first main result we completely settle this conjecture in a far more general form:
\begin{thm}\label{conjecture}
Let $u_k$ be a sequence of blowup solutions of (\ref{dir-0-cond}) with parameter $\lambda_k$ that satisfies (\ref{unif-bound-3}). Then $u_k$ is simple around any blowup point in $\Omega$. 
\end{thm}
Theorem \ref{conjecture} does not impose any symmetry condition on $\Omega$.  As long as $\Omega$ is a bounded open set with smooth boundary, the conclusion of Theorem \ref{conjecture} holds under a natural uniform bound on its mass $ \lambda \int_\Omega e^{u} $. This comes as a complete surprise. (Note that by a nice result of Battaglia \cite{battaglia}, the mass $\lambda \int_\Omega e^u $ is uniformly bounded as long as $\Omega$ is simply connected and $M=1$. In particular when $\Omega=B_1, M=1$, the mass is uniformly bounded. As a result we have completely solved the original conjecture of the first two authors.)  

If the boundary condition is the usual oscillation finiteness assumption (\ref{standard-b}), we can also rule out non-simple blowup in a surprising way. The second main result is:

\begin{thm}\label{main-thm-gen}
Let $\mathfrak{u}_k$ be a sequence of blowup solutions of (\ref{main-e-1}) such that (\ref{standard-a}) and (\ref{standard-b}) hold. If there are at least two blowup points in $\Omega$, each blowup point is a simple blowup point. 
\end{thm}

Theorem \ref{main-thm-gen} is also unexpected. First there is no requirement on the data of $\mathfrak{u}_k$ on $\partial \Omega$ except that the oscillation of $\mathfrak{u}_k$ is finite. There is also no specific requirement of what $\Omega$ has to be. As long as there are at least two blowup points in $\Omega$, each one of them has to be simple. Thus non-simple blow-up, if happens, can only be single.  This greatly simplifies the blowup analysis in many applications.

The conclusions of Theorem \ref{main-thm-gen} and Theorem \ref{conjecture} rule out non-simple blowup phenomenon in several applicable situations. They seem to suggest that the only case that non-simple blowup solutions occur is when the profile of blowup solutions is very close to global solutions in the classification theorem of Prajapat-Tarantello \cite{prajapat}. The proofs of the main results should lead to advances in multiple related problems. Even though we study only one equation in this article, it represents certain situations in systems. For example, in the blowup analysis of Toda systems, which has ties with conformal geometry, algebraic geometry, integrable system and complex analysis \cite{lin-wei-ye,tr-1,tr-2,wu-zhang-siam,zhang-imrn}, one always needs to compare blowup speeds of different components. If different components all tend to infinity in a neighborhood of one blowup point, the behavior of the ``fast" component is similar to a \emph{quantized} singular source to ``slow" components. Therefore the asymptotic behavior of blowup solutions of Liouville equation with quantized singular source provides crucial information for systems.

It is also important to point out that the study of singular equation with ``quantized" singular source is ubiquitous in mathematical literature; the nonsimple blow-up phenomenon also appears in the research of Liouville system \cite{gu-zhang-2}, prescribing Q curvature equation \cite{ahmedou-wu-zhang}, Monge-Ampere equation \cite{chen-lin-cpam-15} and the vortices in a planar model of Euler flows \cite{del-pino-musso}, etc. 

\medskip

The proof of Theorem \ref{main-thm-gen} comes naturally from Theorem \ref{main-thm-1}, which contains most of the key ideas. The proof of Theorem \ref{main-thm-1} is by way of contradiction. If non-simple blowup happens, the blowup solutions (denoted $v_k$) would have $N+1$ local maximums evenly distributed around the unit disk. We shall use $N+1$ global solutions, each is very close to $v_k$ near a local maximum. Even though these $N+1$ global solutions are close to one another, their mutual difference can be captured as kernels of linearized Liouville equation. The mutual locations of the local maximum points plays a key role in our argument. One way to understand this argument is that $v_k$ cannot be very close to different global solutions at the same time. In order to obtain a contradiction, the point-wise estimate has to be very precise. In this article we use Fourier analysis to obtain precise pointwise estimates around one ``easy point" first. Then this accuracy can be passed to other regions by Harnack inequality. In the final step a contradiction can be obtained by comparing Pohozaev identities of $v_k$ and other global solutions. Another key ingredient in our proof is that the behavior of $v_k$ on the boundary of its domain is significantly different from all the approximating global solutions. This difference can be turned into a contrast on coefficient functions, which leads to a contradiction from Pohozaev identities. 
This set of ideas has turned out to be not only successful for single equations \cite{wei-zhang-adv,wei-zhang-plms,wei-zhang-jems}, but also in Toda systems \cite{wei-wu-zhang}. Many related projects such as construction of blowup solutions, the uniqueness of blowup sequence, etc will also be influenced by the proof in this article. The research in these directions will be carried out soon.

The organization of the article is as follows. From Chapter 2 to Chapter 4 we present the key proposition needed for the proof of the main theorem. The proof of Theorem \ref{main-thm-1} in these chapters contains the key ingredients, which are inspired by a series of works of the second and third authors \cite{wei-zhang-adv,wei-zhang-plms,wei-zhang-jems}. Finally we place some computations in the appendix, which consists of the last three chapters. 

\medskip

{\bf Notation:} We will use $B(x_0,r)$ to denote a ball centered at $x_0$ with radius $r$. If $x_0$ is the origin we use $B_r$. $C$ represents a positive constant that may change from place to place.

\section{Local Approximation}
In the section we study the following local equation defined in the unit disk $B_1$ in $\mathbb R^2$:
$$\Delta \mathfrak{u}+e^{\mathfrak{u}}=4\pi N\delta_0, \quad \mbox{in}\quad B_1$$
where $N$ is a positive integer. 
Since $\Delta(\frac 1{2\pi}\log |x|)=\delta_0$ we can use this function to 
write the equation above as
\begin{equation}\label{alter-e}
\Delta u+|x|^{2N}e^{u}=0, 
\end{equation}
if we set $u(x)=\mathfrak{u}(x)-2N\log |x|$. The purpose of this section is to study blowup solutions of (\ref{alter-e}). 

Let $u_k$ be a sequence of solutions of (\ref{alter-e}):
\begin{equation}\label{main-0}
\Delta u_k+|x|^{2N}e^{u_k}=0, \quad \mbox{in}\quad B_1
\end{equation}
 We say $u_k$ is a sequence of blowup solutions with blowup point at the origin, if there exists $x_k\to 0$ such that $u_k(x_k)\to \infty$ as $k\to \infty$.
Suppose the oscillation of $u_k$ on the boundary of $B_1$ is finite:
\begin{equation}\label{BOF}
|u_k(x)-u_k(y)|\le C,\quad \forall x,y\in \partial B_1 
\end{equation}
for some $C>0$ independent of $k$, and there is a uniform bound on the integration of $|x|^{2N}e^{u_k}$:
\begin{equation}\label{unif-energ}\int_{B_1}|x|^{2N}e^{u_k}<C.
\end{equation}

Our goal is to study the asymptotic behavior of $u_k$ near the origin and its relation with the oscillation of $u_k$ on $\partial B_1$ (the boundary of $B_1$). For this purpose we
set
$$\Phi_k(x)=u_k(x)-\frac 1{2\pi}\int_{\partial B_1}u_k, \quad x\in B_1. $$
Since $u_k$ has bounded oscillation on $\partial B_1$, $\Phi_k(0)=0$ and all the derivatives of $\Phi_k$ are uniformly bounded in $B_{1/2}$.  Let $\Phi$ be the limit of $\Phi_k$ over any fixed compact subset of $B_1$. 

Then our assumption of $\Phi_k$ is
\begin{equation}\label{no-oscillation}
\mbox{Either } \Phi\neq 0\quad \mbox{or} \quad \Phi_k\equiv 0.
\end{equation}

\begin{thm}\label{main-thm-1}
Let $u_k$ be a sequence of blowup solutions of (\ref{main-0}) that takes $0$ as its only blowup point in $B_1$. Suppose (\ref{no-oscillation}) and (\ref{unif-energ}) hold. Then $u_k$ is a simple blowup sequence:
$$u_k(x)+2(1+N)\log |x|\le C $$
for some $C>0$. 
\end{thm}

\section{Proof of Theorem \ref{main-thm-1}}
Suppose non-simple blowup does happen. It is well known \cite{kuo-lin-jdg, bart3} that there are exactly $N+1$ local maximums forming a circle around the origin. We use $p_0^k$,...,$p_N^k$  evenly distributed on $\mathbb S^1$ after scaling according to their magnitude: Suppose along a subsequence
$$\lim_{k\to \infty}p_0^k/|p_0^k|=e^{i\theta_0}, $$
then
$$\lim_{k\to \infty} \frac{p_l^k}{|p_0^k|}=e^{i(\theta_0+\frac{2\pi l}{N+1})}, \quad l=1,...,N. $$
For many reasons it is convenient to denote $|p_0^k|$ as $\delta_k$ and define $\bar \mu_k$ as follows:
\begin{equation}\label{muk-dk}
\delta_k=|p_0^k|\quad \mbox{and }\quad \bar \mu_k= u_k(p_0^k)+2(1+N)\log \delta_k.
\end{equation}
 
Since $p_l^k$'s are evenly distributed
around $\partial B_{\delta_k}$, standard results for Liouville equations around a regular blowup point can be applied to have $ u_k(p_l^k)= u_k(p_0^k)+o(1)$. Also,  $\bar \mu_k\to \infty$. The interested readers may look into \cite{kuo-lin-jdg,bart3} for more detailed information.

We write $p_0^k$ as $p_0^k=\delta_ke^{i\theta_k}$ and define $v_k$ as
\begin{equation}\label{v-k-d}
v_k(y)=u_k(\delta_k ye^{i\theta_k})+2(N+1)\log \delta_k,\quad |y|< \delta_k^{-1}.
\end{equation}
If we write out each component, (\ref{v-k-d}) is
$$
v_k(y_1,y_2)=u_k(\delta_k(y_1\cos\theta_k-y_2\sin\theta_k),\delta_k(y_1\sin\theta_k+y_2\cos\theta_k))+2(1+N)\log \delta_k. $$
Then it is standard to verify that $v_k$ solves

\begin{equation}\label{e-f-vk}
\Delta v_k(y)+|y|^{2N}e^{v_k(y)}=0,\quad |y|<\delta_k^{-1},
\end{equation}

Thus the image of $p_0^k$ after scaling is $Q_0^k=e_1=(1,0)$.
Let $Q_1^k$, $Q_2^k$,...,$Q_{N}^k$ be the images of $p_i^k$ $(i=1,...,N)$ after the scaling:
$$Q_l^k=\frac{p_l^k}{\delta_k}e^{-i\theta_k},\quad l=0,...,N. $$
 It is established by Kuo-Lin in \cite{kuo-lin-jdg} and independently by Bartolucci-Tarantello in \cite{bart3} that
\begin{equation}\label{limit-q}
\lim_{k\to \infty} Q_l^k=\lim_{k\to \infty}p_l^k/\delta_k=e^{\frac{2l\pi i}{N+1}},\quad l=0,....,N.
\end{equation}
Then it is proved in our previous work that ( see (3.13) in \cite{wei-zhang-adv})
\begin{equation}\label{Qm-close}
Q_l^k-e^{\frac{2\pi l i}{N+1}}=O(\bar \mu_ke^{-\bar \mu_k}).
\end{equation}
Choosing $3\tau>0$ small and independent of $k$, we can make disks centered at $Q_l^k$ with radius $3\tau$ (denoted as $B(Q_l^k,3\tau ) $) mutually disjoint. The $\bar \mu_k$ in (\ref{muk-dk}) is 
$$\bar \mu_k=\max_{B(Q_0^k,\tau)} v_k.
$$
Since $Q_l^k$ are evenly distributed around $\partial B_1$, it is easy to use standard estimates for single Liouville equations (\cite{zhangcmp,gluck,chenlin1}) to obtain
$$\max_{B(Q_l^k,\epsilon)}v_k=\bar \mu_k+o(1),\quad l=1,...,N. $$

Let
\begin{equation}\label{def-Vk}
V_k(x)=\log \frac{e^{\bar \mu_k}}{(1+\frac{e^{\bar \mu_k}}{8(1+N)^2}|y^{N+1}-e_1|^2)^2}.
\end{equation}
Clearly $V_k$ is a solution of
\begin{equation}\label{eq-for-Vk}
\Delta V_k+|y|^{2N}e^{V_k}=0,\quad \mbox{in}\quad \mathbb R^2, \quad V_k(e_1)=\mu_k.
\end{equation}
This expression is based on the classification theorem of Prajapat-Tarantello \cite{prajapat}.  Now we use the following expansion of $V_k$ for $|y|=L_k$ ($L_k=\delta_k^{-1}$)
\begin{align}\label{expan-Vk-Lk}
V_k(y)&=-\bar \mu_k+2\log (8(N+1)^2)-4(N+1)\log L_k+\frac{2}{L_k^{2N+2}}\\
&+\frac{4\cos((N+1)\theta)}{L_k^{N+1}}
+\frac{4}{L_k^{2N+2}}\cos((2N+2)\theta)\nonumber\\
&+O(L_k^{-3N-3})+O(e^{-\bar \mu_k}L_k^{-2N-2}). \nonumber
\end{align}

To eliminate the main oscillation of $v_k-V_k$ on $\partial \Omega_k$, we set
$\phi_{v,k}(\delta_k\cdot )$ be the harmonic function on $\partial \Omega_k$ such that $V_k(y)-\phi_{v,k}(\delta_ky)$ has no oscillation. The reason we use $\phi_{v,k}(\delta_k\cdot)$ is because $\phi_{v,k}$ is the harmonic function that eliminates the oscillation of $\tilde V_k$ on $\partial B_1$. $\tilde V_k$ is a re-scaled $V_k$: 
$$\tilde V_k(\delta_k y)+2(1+N)\log \delta_k=V_k(y). $$
on $\partial B_1$. From the expression of $V_k$ we see that $\phi_{v,k}\to 0$ uniformly on $B_1$. It is easy to see that the leading terms of $\phi_{v,k}(\delta_ky)$ are
$$\phi_{v,k}(\delta_ky)=\frac{4}{L_k^{2N+2}}r^{N+1}\cos((N+1)\theta)+\frac{4}{L_k^{4N+4}}r^{2N+2}\cos((2N+2)\theta)+.... $$

Naturally we define $\phi_{0,k}(y)$, $h_{0,k}$ and $v_{0,k}$ as
\begin{align}\label{har-osci}&\phi_{0,k}(y)=\Phi_k(\delta_ky)-\phi_{v,k}(\delta_ky)\\
&=\Phi_k(\delta_ky)-4\delta_k^{2N+2}r^{N+1}\cos((N+1)\theta)
-\delta_k^{4N+4}r^{2N+2}\cos((2N+2)\theta)+...\nonumber
\end{align}
\begin{equation}\label{h0k-def}
h_{0,k}=e^{\phi_{0,k}}. 
\end{equation} 
$$v_{0,k}=v_k-\phi_{0,k},$$
and we write the equation of $v_{0,k}$ as 
\begin{equation}\label{e-v0k}\Delta v_{0,k}+h_{0,k}|y|^{2N}e^{v_{0,k}}=0, \quad \mbox{in}\quad \Omega_k 
\end{equation}

Based on the definition of $h_{0,k}$ in (\ref{h0k-def}) we prove
\begin{lem}\label{key-lem-2} There exist and integer $L>0$ ,$\delta_k^*\in (\delta_k^L,\delta_k)$ an integer $0\le s\le N$ such that
\begin{equation}\label{non-vanish}
\nabla h_{0,k}(e^{\frac{2\pi is}{N+1}})/\delta_k^{*}\neq 0.
\end{equation}
\end{lem}

{\bf Proof of Lemma \ref{key-lem-2}:}
The proof of (\ref{non-vanish}) is based on (\ref{no-oscillation}).
We first write the expansion of $\Phi_k(x)$ as
$$\Phi_k(x)=\sum_{n=1}^{\infty}r^n(a_n^k\cos(n\theta)+b_n^k\sin(n\theta))$$
where $x=r e^{i\theta}$ and the two sequences $|a_n^k|$ and $|b_n^k|$ are uniformly bounded from above. They could tend to zero for some $n$. The expansion of $\phi_{v,k}(x)$ can be written as
$$\phi_{v,k}(x)=\sum_{n=1}^{\infty}r^n(a_{n,v}^k\cos(n\theta)+b_{n,v}^k\sin(n\theta)) $$
with $a_{n,v}^k$ and $b_{n,v}^k$ both tending to zero as $k\to \infty$ for each fixed $n$. 
Then clearly 
\begin{align*}
   \phi_{0,k}(y)&=\Phi_k(\delta_ky)-\phi_{v,k}(\delta_ky)\\
   &=\sum_{n=1}^{\infty}\delta_k^n r^n\bigg ((a_{n}^k-a_{n,v}^k)\cos n\theta+
(b_n^k-b_{n,v}^k)\sin n\theta)\bigg ). 
\end{align*}
Obviously we use $y=r e^{i\theta}$ in this case. 
If $\Phi\neq 0$ (see \ref{no-oscillation}) 
there exists an integer $L$ such that $|a_L^k|+|b_L^k|\ge C>0$ for some $C>0$ independent of $k$. Now we set $\delta_k^*$ as 
$$\delta_k^*:=\sum_{n=1}^L\delta_k^n(|a_n^k-a_{n,v}^k|+|b_n^k-b_{n,v}^k|). $$
We make three trivial remarks about $\delta_k^*$. First 
$$\delta_k^*\ge C\delta_k^L $$
because $|a_L^n|+|b_L^n|\ge C$ for some $C$ independent of $k$. Second, since it is only a finite sum in the definition of $\delta_k^*$, $\delta_k^*$ is comparable to 
$$\max_{1\le n\le L}(|a_n^k-a_{n,v}^k|+|b_n^k-b_{n,v}^k|)\delta_k^n.$$
The third observation is 
$$\sum_{n=L+1}^{\infty}(|a_n^k-a_{n,v}^k|+|b_n^k-b_{n,v}^k|)\delta_k^n=O(\delta_k^{L+1}). $$
In other words, the quantity is far less than $\delta_k^*$. 
Our goal is to compute $\nabla \phi_{0,k}$ at points close to $\partial B_1$. We use the following formula:
$$|\nabla \phi_{0,k}(x)|^2=|\partial_r\phi_{0,k}|^2+\frac 1{r^2}|\partial_{\theta}\phi_{0,k}|^2, \quad |x|\sim 1. $$
For convenience we use 
$$(a_n^k-a_{n,v}^k)\delta_k^n=\delta_k^*|c_n^k|\cos(\beta_n^k), \quad (b_n-b_{n,v})\delta_k^n=\delta_k^*|c_n^k|\sin(\beta_n^k). $$
Note that at least one $|c_n^k|\sim 1$. In fact we consider those comparable to $1$ in the summation from $1$ to $L$. Then
$$\partial_r\phi_{0,k}=\delta_k^*\sum_{n=1}^L n r^{n-1}|c_n^k|\cos(\beta_n^k-n\theta)+O(\delta_k^{L+1})$$
\begin{align*}
 &\partial_{\theta}\phi_{0,k}\\
 =&\delta_k^*\sum_{n=1}^Lr^n(-n|c_n^k|\cos(\beta_n^k)\sin(n\theta)+n|c_n^k|\sin(\beta_n^k)\cos(n\theta))+O(\delta_k^{L+1})\\
 =&\delta_k^*\sum_{n=1}^N n r^n |c_n^k| \sin(\beta_n^k-n\theta) +O(\delta_k^{L+1}).
\end{align*}
\begin{align}\label{grad-Phi}
    &|\nabla \phi_{0,k}|^2
    =|\partial_r\phi_{0,k}|^2+\frac{1}{r^2}|\partial_{\theta}\phi_{0,k}|^2\\
    =&(\delta_k^*)^2(\sum_{n=1}^Ln^2r^{2n-2}|c_n^k|^2+2\sum_{s<t}st r^{s+t-2}|c_s^k||c_t^k|\cos(\beta_s^k-\beta_t^k-(s-t)\theta)\nonumber\\
    &+O(\delta_k^{L+1}).\nonumber
\end{align}
Now we have the freedom to choose $e^{\frac{2\pi l i}{N+1}}$. Since all these points are on the unit circle, $r=1$ in the evaluation of (\ref{grad-Phi}). By Cauchy's inequality we see that if $|\nabla \phi_{0,k}|=O(\delta_k^{L+1})$ at one of $e^{\frac{2\pi i l}{N+1}}$ it is obviously not the case on another such point, which means at that point, say $e^{\frac{2\pi i l}{N+1}}$, $|\nabla \phi_{0,k}(e^{\frac{2\pi i l}{N+1}})|$ is comparable to $\delta_k^*$.  Lemma \ref{key-lem-2} is established if $\Phi\neq 0$.

Finally if $\Phi_k\equiv 0$, based on the expansion of $V_k$ we take $\delta_k^*=\delta_k^{2N+2}$. 
Lemma \ref{key-lem-2} is established. $\Box$

\medskip

If we let $Q_0^k$ denote the local maximum point of $v_{0,k}$, the difference between $Q_0^k$ and $e_1$ is $O(\delta_k^*e^{-\bar\mu_k})$ by the non-degeneracy of $v_k$ around $e_1$ (Indeed, if we use $\bar \epsilon_k=e^{-\bar \mu_k/2}$ to be the scaling factor, after scaling, it is easy to see that the location of the critical point is $O(\bar \epsilon_k \delta_k^*)$ for the scaled function. In the original setting before the scaling, the location of the local maximum is $O(\delta_k^*e^{-\bar \mu_k})$ away from $e_1$. The maximum of $v_{0,k}$ is $\bar \mu_k+O(\delta_k^*)$. 

Let $V_{0,k}$ be the global solution of that agrees with $v_{0,k}$ at $Q_0^k$:
\begin{equation}\label{sta-1}
v_{0,k}(Q_0^k)=V_{0,k}(Q_0^k),\quad \nabla v_{0,k}(Q_0^k)=\nabla V_{0,k}(Q_0^k)=0
\end{equation}
and the equation for $V_{0,k}$ is 
\begin{equation}\label{V0k-eq}
\Delta V_{0,k}+|y|^{2N}h_{0,k}(Q_0^k)e^{V_{0,k}}=0,\quad \mbox{in}\quad \mathbb R^2,\quad \int_{\mathbb R^2} |y|^{2N}e^{V_{0,k}}<\infty. 
\end{equation}
Of course $V_{0,k}$ is a small perturbation of the previous $V_k$, and we see that the oscillation of $v_{0,k}-V_{0,k}$ is $O(\delta_k^*\delta_k^{N+1}e^{-\bar \mu_k})$ on $\partial \Omega_k$. At this moment we set the scaling factor to be
$$\epsilon_k=e^{-\frac 12 \mu_k}. \quad \mu_k=v_{0,k}(Q_0^k) $$
Indeed, the expression of $V_{0,k}$ is 
$$V_{0,k}(y)=\log \frac{e^{\mu_k}}{(1+\frac{h_{0,k}(Q_0^k) e^{\mu_k}}{D}|y^{N+1}-1-p_0^k|^2)^2} $$
where $p_0^k=O(e^{-\mu_k}\delta_k^*)$. 
From here we see that the perturbation of $v_{0,k}-V_{0,k}$ on $\partial \Omega_k$ is $O(\delta_k^*\delta_k^{N+1}e^{-\mu_k})$:
\begin{equation}\label{bry-small}
|(v_{0,k}-V_{0,k})(x)-(v_{0,k}-V_{0,k})(y)|\le C\delta_k^*\delta_k^{N+1}e^{-\mu_k},\quad \forall x,y\in\partial \Omega_k. 
\end{equation}

\subsection{Point-wise estimate for $v_{0,k}-V_{0,k}$}

In the appendix we establish the closeness of local maximum points of $v_k$ with $e^{\frac{2i\pi l}{N+1}}$ for $l=0,...,N$. As mentioned before the local maximums of $v_{0,k}$ is only $O(\delta_k^*\epsilon_k^2)$ perturbation of the corresponding local maximum points of $v_k$. 

Another observation is that based on (\ref{Qm-close}) we have 
$$\epsilon_k^{-1} | Q_l^k-e^{i\beta_l}|\le C\epsilon_k^{\epsilon},\quad l=0,...,N$$
for some small $\epsilon>0$. Thus $\xi_k$ tends to $U$ after scaling. We need this fact in our argument.

Now we cite Proposition 3.1 of \cite{wei-zhang-plms}: 

\emph{Proposition 3.1 of \cite{wei-zhang-plms}: Let $l=0,...,N$ and $\delta$ be small so that $B(e^{i\beta_l},\delta)\cap B(e^{i\beta_s},\delta)=\emptyset$ for $l\neq s$.
In each $B(e^{i\beta_l},\delta)$
\begin{equation}\label{global-close}
|v_{0,k}(x)-V_{0,k}(x)|\le \left\{\begin{array}{ll}
C\mu_ke^{-\mu_k/2},\quad |x-e^{i\beta_l}|\le Ce^{-\mu_k/2}, \\
\\
C\frac{\mu_ke^{-\mu_k}}{|x-e^{i\beta_l}|}+O(\mu_k^2e^{-\mu_k}),\quad Ce^{-\mu_k/2}\le |x-e^{i\beta_l}|\le \delta.
\end{array}
\right.
\end{equation} }

\medskip

\begin{rem} We only need a re-scaled version of Proposition 3.1 of \cite{wei-zhang-plms}:  
\begin{equation}\label{vk-Vk-2}
|v_{0,k}(e^{i\beta_l}+\epsilon_ky)-V_{0,k}(e^{i\beta_l}+\epsilon_ky)|\le C\epsilon_k^{\epsilon} (1+|y|)^{-1},\quad 0<|y|<\tau \epsilon_k^{-1}.
\end{equation}
for some small constants $\epsilon>0$ and $\tau>0$ both independent of $k$. 
\end{rem}

\begin{rem} The main idea of the proof of Theorem \ref{main-thm-1} can be observed from the equations of $v_{0,k}$ (\ref{e-v0k}) and $V_{0,k}$ (\ref{V0k-eq}). While (\ref{V0k-eq}) has a constant coefficient $h_{0,k}(Q_0^k)$, (\ref{e-v0k}) has a function $h_{0,k}$ whose derivative is not zero at some $Q_s^k$. We shall obtain a precise pointwise estimate of $v_{0,k}-V_{0,k}$ and the difference on these coefficient functions will lead to a contradiction. 
\end{rem}

\begin{prop}\label{key-w8-8} Let $w_{0,k}=v_{0,k}- V_{0,k}$, then
$$|w_{0,k}(y)|\le C\delta_k^* \quad y\in \Omega_k:=B(0,\delta_k^{-1}), $$
\end{prop}

\noindent{\bf Proof of Proposition \ref{key-w8-8}:}
Let $M_k$ be the maximum of $|w_{0,k}|$ on $\Omega_k$. 
By way of contradiction we assume that
\begin{equation}\label{big-mk}
M_k/(\delta_k^*)\to \infty. 
\end{equation}

Now we recall the equation for $v_{0,k}$ is (\ref{e-f-vk}),
$v_{0,k}-V_{0,k}$ is very close to a constant on $\partial B(0,\delta_k^{-1})$. Moreover
\begin{equation}\label{control-e}
w_{0,k}(Q_0^k)= |\nabla w_{0,k}(Q_0^k)|=0. 
\end{equation}
Recall that $V_{0,k}$ defined in (\ref{def-Vk}) satisfies (\ref{V0k-eq}).

We shall derive a precise, point-wise estimate of $w_{0,k}$ in $B_3\setminus \cup_{l=1}^{N}B(Q_l^k,\tau)$ where $\tau>0$ is a small number independent of $k$. We shall prove that $w_{0,k}$ is very small in $B_3$ if we exclude all bubbling disks except the one around $e_1$. 

Now based on (\ref{e-v0k}) and (\ref{V0k-eq}) we write the equation of $w_{0,k}$ as
\begin{equation}\label{eq-wk}
\Delta w_{0,k}+|y|^{2N}h_{0,k}(y)e^{\xi_k}w_{0,k}=|y|^{2N}(h_{0,k}(Q_0^k)-h_{0,k}(y))e^{V_{0,k}}
\end{equation}
 in $\Omega_k$, where $\xi_k$ is obtained from the mean value theorem:
$$
e^{\xi_k(x)}=\left\{\begin{array}{ll}
\frac{e^{v_{0,k}(x)}-e^{ V_{0,k}(x)}}{v_{0,k}(x)- V_{0,k}(x)},\quad \mbox{if}\quad v_{0,k}(x)\neq V_{0,k}(x),\\
\\
e^{ V_{0,k}(x)},\quad \mbox{if}\quad v_{0,k}(x)=V_{0,k}(x).
\end{array}
\right.
$$
An equivalent form is
\begin{equation}\label{xi-k}
e^{\xi_k(x)}=\int_0^1\frac d{dt}e^{tv_{0,k}(x)+(1-t)V_{0,k}(x)}dt=e^{ V_{0,k}(x)}\Big(1+\frac 12w_{0,k}(x)+O(w_{0,k}(x)^2)\Big).
\end{equation}
Note that the oscillation of $w_{0,k}$ on $\partial \Omega_k$ is $O(\delta_k^*\delta_k^{N+1}e^{-\mu_k})$. A trivial observation is that the right hand side of (\ref{eq-wk}) is zero if $y=Q_0^k$. This simple fact determines a pathway for our argument. We will prove smallness of $w_k$ around $Q_0^k$ first and then pass it to other places. 

By normalizing $w_{0,k}$ we shall study the following function: 
$$\tilde w_k(y)=w_{0,k}(y)/M_k,\quad x\in \Omega_k. $$

Clearly $\max_{x\in \bar\Omega_k}|\tilde w_k(x)|=1$. The equation for $\tilde w_k$ is
\begin{equation}\label{t-wk}
\Delta \tilde w_k(y)+|y|^{2N}h_{0,k}(y)e^{\xi_k}\tilde w_k(y)=O(\sigma_k)|y|^{2N}(y-Q_0^k)e^{V_{0,k}}+O(\sigma_k)|y-Q_0^k|^2e^{V_{0,k}}
\end{equation}
in $\Omega_k$. Here $\delta_k^*/M_k=\sigma_k\to 0$. 
Also on the boundary, the oscillation of $\tilde w_k$ is $o(\delta_k^{N+1}e^{-\mu_k})$.
By Proposition 3.1 of \cite{wei-zhang-plms}
\begin{equation}\label{xi-V-c}
\xi_k(Q_0^k+\epsilon_k z)=V_{0,k}(Q_0^k+\epsilon_kz)+O(\epsilon_k^{\epsilon})(1+|z|)^{-1}
\end{equation}

Since $V_{0,k}$ is not exactly symmetric around $Q_0^k$, we shall replace the re-scaled version of $V_{0,k}$ around $Q_0^k$ by a radial function.
Let $U_k$ be solutions of
\begin{equation}\label{global-to-use}
\Delta U_k+h_{0,k}(Q_0^k)e^{U_k}=0,\quad \mbox{in}\quad \mathbb R^2, \quad U_k(0)=\max_{\mathbb R^2}U_k=0.
\end{equation}
By the classification theorem of Caffarelli-Gidas-Spruck \cite{CGS} we have
$$U_k(z)=\log \frac{1}{(1+\frac{h_{0,k}(Q_0^k)}{8}|z|^2)^2}$$
and standard refined estimates yield (see \cite{chenlin1,zhangcmp,gluck})
\begin{equation}\label{Vk-rad}
V_{0,k}(Q_0^k+\epsilon_k z)+2\log \epsilon_k=U_k(z)+O(\epsilon_k)|z|+O(\mu_k^2\epsilon_k^2).
\end{equation}
Also we observe that
\begin{equation}\label{log-rad}
\log |Q_0^k+\epsilon_k z|=O(\epsilon_k)(1+|z|).
\end{equation}

Thus, the combination of (\ref{xi-V-c}), (\ref{Vk-rad}) and (\ref{log-rad}) gives
\begin{align}\label{xi-U}
&2N\log |Q_0^k+\epsilon_kz|+\xi_k(Q_0^k+\epsilon_k z)+2\log \epsilon_k-U_k(z)\\
=&O(\epsilon_k^{\epsilon})(1+|z|)\quad 0\le |z|<\tau \epsilon_k^{-1}.
 \nonumber
\end{align}
for a small $\epsilon>0$ independent of $k$. 
Since we shall use the re-scaled version, based on (\ref{xi-U}) we have
\begin{equation}\label{xi-eU}
\epsilon_k^2 |Q_0^k+\epsilon_k z|^{2N}h_{0,k}(Q_0^k+\epsilon_k z)e^{\xi_k(Q_0^k+\epsilon_k z)}
=h_{0,k}(Q_0^k) e^{U_k(z)}+O( \epsilon_k^{\epsilon})(1+|z|)^{-3}
\end{equation}
Here we note that the estimate in (\ref{xi-U}) is not optimal.  In the following we shall put the proof of Proposition \ref{key-w8-8} into a few estimates. In the first estimate we prove

\begin{lem}\label{t-w-1-better} 
\begin{equation}\label{for-lambda-k}
|\tilde w_k(Q_0^k+\epsilon_kz)|\le o(\epsilon_k) (1+|z|),\quad 0<|z|<\tau \epsilon_k^{-1}.
\end{equation}
for some $\tau>0$.
\end{lem}

\noindent{Proof of Lemma \ref{t-w-1-better}:}
{\bf Step one:} In this step we prove the following statement:
For $\delta>0$ small and independent of $k$,
\begin{equation}\label{key-step-1}
\tilde w_k(y)=o(1),\quad \nabla \tilde w_k=o(1) \quad \mbox{in}\quad B(e_1,\delta)\setminus B(e_1,\delta/8)
\end{equation}
where $B(e_1,3\delta)$ does not include other blowup points.

If (\ref{key-step-1}) is not true, we have, without loss of generality that $\tilde w_k\to c\neq 0$. This is based on the fact that $\tilde w_k$ tends to a global harmonic function with removable singularity. So $\tilde w_k$ tends to constant. Let
\begin{equation}\label{w-ar-e1}
W_k(z)=\tilde w_k(Q_0^k+\epsilon_kz), \quad \epsilon_k=e^{-\frac 12 v_{0,k}(Q_0^k)},
\end{equation}
then if we use $W$ to denote the limit of $W_k$, we have
$$\Delta W+e^UW=0, \quad \mathbb R^2, \quad |W|\le 1, $$
and $U$ is a solution of $\Delta U+e^U=0$ in $\mathbb R^2$ with $\int_{\mathbb R^2}e^U<\infty$. Since $0$ is the local maximum of $U$,
$$U(z)=\log \frac{1}{(1+\frac 18|z|^2)^2}. $$
Here we further claim that $W\equiv 0$ in $\mathbb R^2$ because $W(0)=|\nabla W(0)|=0$, a fact well known based on the classification of the kernel of the linearized operator. Going back to $W_k$, we have
$$W_k(z)=o(1),\quad |z|\le R_k \mbox{ for some } \quad R_k\to \infty. $$

Based on the expression of $\tilde w_k$, (\ref{Vk-rad}) and (\ref{xi-eU}) we write the equation of $W_k$ as
\begin{equation}\label{e-Wk}
\Delta W_k(z)+h_{0,k}(Q_0^k)e^{U_k(z)}W_k(z)=E_k,
\end{equation}
for $|z|<\delta_0 \epsilon_k^{-1}$ where a crude estimate of the error term $E_k$ is
\begin{equation*}
E_k(z)=o(1)\epsilon_k^{\epsilon}(1+|z|)^{-3}.
\end{equation*}

Let
\begin{equation}\label{for-g0}
g_0^k(r)=\frac 1{2\pi}\int_0^{2\pi}W_k(r,\theta)d\theta.
\end{equation}
Then clearly $g_0^k(r)\to c>0$ for $r\sim \epsilon_k^{-1}$.
 The equation for $g_0^k$ is
\begin{align*}
&\frac{d^2}{dr^2}g_0^k(r)+\frac 1r \frac{d}{dr}g_0^k(r)+h_{0,k}(Q_0^k)e^{U_k(r)}g_0^k(r)= E_0^k(r)\\
&g_0^k(0)=\frac{d}{dr}g_0^k(0)=0.
\end{align*}
where $E_0^k(r)$ has the same upper bound as that of $E_k(r)$:
$$|E_0^k(r)|\le o(1)\epsilon_k^{\epsilon}(1+r)^{-3}. $$

For the homogeneous equation, the two fundamental solutions are known: $g_{01}^k$, $g_{02}^k$, where
$$g_{01}^k=\frac{1-c_kr^2}{1+c_kr^2},\quad c_k=\frac{h_{0,k}(Q_0^k)}8.$$
By the standard reduction of order process, $g_{02}^k(r)=O(\log r)$ for $r>1$ with a bound independent of $k$.
Then it is easy to obtain, assuming $|W_k(z)|\le 1$, that
\begin{align*}
|g_0^k(r)|\le C|g_{01}^k(r)|\int_0^r s|E_0^k(s) g_{02}^k(s)|ds+C|g_{02}^k(r)|\int_0^r s|g_{01}^k(s) E_0^k(s)|ds\\
\le C\epsilon_k^{\epsilon}\log (2+r). \quad 0<r<\delta_0 \epsilon_k^{-1}.
\end{align*}
Clearly this is a contradiction to (\ref{for-g0}). We have proved $c=0$, which means $\tilde w_k=o(1)$ in $B(e_1, \delta_0)\setminus B(e_1, \delta_0/8)$.
Then it is easy to use the equation for $\tilde w_k$ and standard Harnack inequality to prove
$\nabla \tilde w_k=o(1)$ in the same region.
(\ref{key-step-1}) is established. 

\medskip

{\bf Step two}: Now we extend the estimate to the whole neighborhood of $e_1$. The estimate will be obtained in a progressive way. Let $W_k$ be defined as in (\ref{w-ar-e1}). In order to obtain a better estimate we need to write the equation of $W_k$ more precisely than (\ref{e-Wk}):
\begin{equation}\label{w-more}
\Delta W_k+h_{0,k}(Q_0^k)e^{\Theta_k}W_k=o(\sigma_k)\epsilon_ky(1+|y|)^{-4}, \quad z\in B(0,\tau \epsilon_k^{-1})
\end{equation}
where
$\Theta_k$ is defined by
$$e^{\Theta_k(z)}=|Q_0^k+\epsilon_k z|^{2N}e^{\xi_k(Q_0^k+\epsilon_kz)+2\log \epsilon_k}. $$

Here we observe that by step one $W_k=o(1)$ on $\partial B(0,\tau \epsilon_k^{-1})$. 
In the computation of (\ref{w-more}) we also used 
$$h_{0,k}(Q_0^k+\epsilon_ky)=h_{0,k}(Q_0^k)+O(\epsilon_k y). $$
To replace $\Theta_k$ by $U_k$ we have an extra error that depends on the bound of $W_k$:
\begin{equation}\label{w-more-2}
\Delta W_k+h_{0,k}(Q_0^k)e^{U_k}W_k=E_k, \quad z\in B(0,\tau\epsilon_k^{-1})
\end{equation}
where $E_k=E_1+E_2$ and $E_1$ is the right hand side of (\ref{w-more}). Thus 
$$|E_1|\le C\sigma_k\epsilon_ky(1+|y|)^{-4}\quad |E_2|\le C\sigma_k\epsilon_k^{\epsilon}(1+|z|)^{-3}. $$
Note that in the crude bound of $E_2$, we used $|W_k|\le C$. If the bound of $|W_k|$ is better, the estimate of $E_2$ will improve as well. 

The proof  is by considering the projection of $W_k$ in different nodes. Let $g_0^k$ be the projection on $1$, then we have
\begin{align*}
&\frac{d^2}{dr^2}g_0^k(r)+\frac 1r \frac{d}{dr}g_0^k(r)+h_{0,k}(Q_0^k)e^{U_k}g_0^k=E_{k,0}\\
&g_0^k(0)=\frac{d}{dr}g_0^k(0)=0
\end{align*}
where $E_{k,0}$ is the projection of $E_k$ to $1$. Just like in the proof of step one, the estimate of $g_0^k$ is
\begin{align*}
|g_0^k(r)|&\le C|g_{01}(r)|\int_0^r  s|E_{k,0}(s)||g_{02}^k(s)|ds\\
&+C|g_{02}^k(r)|\int_0^r s|E_{k,0}(s)||g_{01}^k(s)|ds
\end{align*}
where
$$g_{01}^k(r)=\frac{1-c_k r^2}{1+c_k r^2}, \quad c_k=\frac{h_{0k}(Q_0^k)}8. $$
and $g_{02}^k=O(\log r)$ for $r>1$ and for $r$ close to $0$. If we just use 
$$|E_{k,0}(s)|\le o(1)\epsilon_k^{\epsilon}(1+s)^{-3}, $$
we have
$$|g_0^k(r)|=o(1)\epsilon_k^{\epsilon}\log (2+r). $$
Let $g_1^k$ be the radial part of the projection of $W_k$ on $e^{i\theta}$. 
Then $g_1^k$ satisfies
$$\frac{d^2}{dr^2}g_1^k(r)+\frac 1r \frac{d}{dr}g_1^k(r)+(h_{0,k}(Q_0^k)e^{U_k}-\frac 1{r^2})g_1^k(r)=E_{k,1}$$
and we use the following crude upper bound for $E_{k,1}$
$$E_{k,1}(r)=o(1)\epsilon_k^{\epsilon}(1+r)^{-3} $$
and $g_1^k(0)=\frac{d}{dr}g_1^k(0)=0$. Then the two fundamental solutions are
$$g_{11}^k(r)=\frac{r}{1+c_k r^2}, $$
$g_{12}^k(r)$ behaves like $O(r^{-1})$ near $0$ and behaves like $O(r)$ near infinity. 
\begin{align*}
|g_1^k(r)|\le &o(1)|g_{11}^k(r)|\int_0^r s \epsilon_k^{\epsilon}(1+s)^{-3}|g_{12}^k(s)| ds\\
&+o(1)|g_{12}^k(r)|\int_0^r s 
\epsilon_k^{\epsilon} (1+s)^{-3} |g_{11}^k(s)|ds\\
&=o(1)\epsilon_k^{\epsilon}(1+r)
\end{align*} 

For $l\ge 2$, let $g_l^k$ be the radial part of the projection on $e^{il\theta}$, then the equation  is
$$\frac{d}{dr^2}g_l^k(r)+\frac 1r\frac{d}{dr}g_l^k(r)+(h_{0,k}(Q_0^k)e^{U_k}-\frac{l^2}{r^2})g_l^k(r)=o(1)\epsilon_k^{\epsilon}(1+r)^{-3} $$
with $g_l^k(0)=0$ and $g_l^k(L_k)=s_{l,k}$ for some $s_{l,k}=o(1)$ as a result of step one. 
For the corresponding homogeneous equation, the two fundamental solutions $g_{l1}^k$, $g_{l2}^k$ can be chosen to behave like
$$g_{l1}^k\sim r^l \mbox{near $0$ and $\infty$}, $$
$$g_{l2}^k\sim r^{-l} \mbox{near $0$ and $\infty$} $$
with bounds independent of $k$. 
Then 
\begin{align*}
 | g_l^k(r)|\le |c_{1,l}^k g_{l1}^k(r)|+o(\epsilon_k^{\epsilon})|g_{l1}^k(r)|\int_r^{\infty}\frac{s}{l}(1+s)^{-3}|g_{l2}^k(s)|ds\\
  +o(\epsilon_k^{\epsilon})|g_{l2}^k(r)|\int_0^r\frac{s}{l}|g_{l1}^k(s)|(1+s)^{-3}ds  \\
  \le |c_{l,1}^kg_{l1}^k(r)|+\frac{o(\epsilon_k^{\epsilon})}{l^2}.
\end{align*}
By $g_{l}(\tau \epsilon_k^{-1})=s_{l,k}$ we have 
$|c_{l,1}^k|\le C s_{l,k}\epsilon_k^l$. Thus the summation of projections on all nodes is convergent and the summation leads to 
$$|W_k(z)|\le o(1)\epsilon_k^{\epsilon}(1+|z|). $$
Using this new bound in the previous machinery we will obtain 
$$|W_k(y)|\le o(1)\epsilon_k^{2\epsilon}(1+|z|). $$
A repetitive application of this process eventually makes $E_1$ the leading error term in $E_k$  and we have
$$|W_k(z)|\le o(1)\epsilon_k(1+|z|),\quad |z|\le \tau \epsilon_k^{-1}. $$
Lemma \ref{t-w-1-better} is
established. $\Box$

\medskip

The smallness of $\tilde w_k$ around $e_1$ can be used to obtain the following key estimate:
\begin{lem}\label{small-other}
\begin{equation}\label{key-step-2}
\tilde w_k=o(1)\quad \mbox{in}\quad B(e^{i\beta_l},\tau)\quad l=1,..,N.
\end{equation}
\end{lem}

\noindent{\bf Proof of Lemma \ref{small-other}:}
We abuse the notation $W_k$ by defining it as
$$W_k(z)=\tilde w_k(e^{i\beta_l}+\epsilon_k z),\quad z\in B(0,\tau \epsilon_k^{-1}). $$
Here we point out that based on (\ref{Qm-close}) we have $\epsilon_k^{-1}|Q_l^k-e^{i\beta_l}|\to 0$. So the scaling around $e^{i\beta_l}$ or $Q_l^k$ does not affect the
limit function.

$$\epsilon_k^2 |e^{i\beta_l}+\epsilon_kz|^{2N}h_{0,k}(e^{i\beta_l}+\epsilon_kz)e^{\xi_k(e^{i\beta_l}+\epsilon_kz)}\to e^{U(z)} $$
where $U(z)$ is a solution of
$$\Delta U+e^U=0,\quad \mbox{in}\quad \mathbb R^2, \quad \int_{\mathbb R^2}e^U<\infty. $$
Since $W_k$ converges to a solution of the linearized equation:
$$\Delta W+e^UW=0, \quad \mbox{in}\quad \mathbb R^2. $$
$W$ can be written as a linear combination of three functions:
$$W(x)=c_0\phi_0+c_1\phi_1+c_2\phi_2, $$
where
$$\phi_0=\frac{1-\frac 18 |x|^2}{1+\frac 18 |x|^2} $$
$$\phi_1=\frac{x_1}{1+\frac 18 |x|^2},\quad \phi_2=\frac{x_2}{1+\frac 18|x|^2}. $$
The remaining part of the proof consists of proving $c_0=0$ and $c_1=c_2=0$. 

\noindent{\bf Step one: $c_0=0$.}
First we write the equation for $W_k$ in a convenient form. Since
$$|e^{i\beta_l}+\epsilon_kz|^{2N}=1+O(\epsilon_k z),$$
and
$$\epsilon_k^2h_{0,k}(e^{i\beta_l}+\epsilon_kz)e^{\xi_k(e^{i\beta_l}+\epsilon_kz)}=e^{U_k(z)}+O(\epsilon_k^{\epsilon})(1+|z|)^{-3}. $$
Based on (\ref{t-wk}) we write the equation for $W_k$ as
\begin{equation}\label{around-l}
\Delta W_k(z)+e^{U_k}W_k=E_l^k(z)
\end{equation}
where
$$E_l^k(z)=O(\epsilon_k^{\epsilon})(1+|z|)^{-3}\quad \mbox{in}\quad \Omega_{k,l}.$$
In order to prove $c_0=0$, the key is to control the derivative of $W_0^k(r)$ where
$$W_0^k(r)=\frac 1{2\pi r}\int_{\partial B_r} W_k(re^{i\theta})dS, \quad 0<r<\tau \epsilon_k^{-1}. $$
To obtain a control of $\frac{d}{dr}W_0^k(r)$ we use $\phi_0^k(r)$ as the radial solution of
\begin{equation}\label{radial-1}\Delta \phi_0^k+e^{U_k}\phi_0^k=0, \quad \mbox{in }\quad \mathbb R^2. \end{equation}
When $k\to \infty$, $\phi_0^k\to \phi_0$. From multiplying $\phi_0^k$ to (\ref{around-l}) and multiplying $W_k$ to (\ref{radial-1}) we have
\begin{equation}\label{c-0-pf}
\int_{\partial B_r}(\partial_{\nu}W_k\phi_0^k-\partial_{\nu}\phi_0^kW_k)=o(\epsilon_k^{\epsilon}). \end{equation}

Thus from (\ref{c-0-pf}) we have
\begin{equation}\label{W-0-d}
\frac{d}{dr}W_0^k(r)=\frac{1}{2\pi r}\int_{\partial B_r}\partial_{\nu}W_k=o(\epsilon_k^{\epsilon})/r+O(1/r^3),\quad 1<r<\tau \epsilon_k^{-1}.
\end{equation}
Since we have known that
$$W_0^k(\tau \epsilon_k^{-1})=o(1). $$
By the fundamental theorem of calculus we have
$$W_0^k(r)=W_0^k(\tau\epsilon_k^{-1})+\int_{\tau \epsilon_k^{-1}}^r\Big(\frac{o(\epsilon_k^{\epsilon})}{s}+O(s^{-3})\Big)ds=O(1/r^2)+O(\epsilon_k^{\epsilon}\log
\frac{1}{\epsilon_k}) $$
for $r\ge 1$. Thus
$c_0=0$ because $W_0^k(r)\to c_0\phi_0$, which means when $r$ is large, it is $-c_0+O(1/r^2)$.

\medskip

\noindent{\bf Step two $c_1=c_2=0$}.
 We first observe that Lemma \ref{small-other} follows from this. Indeed, once we have proved $c_1=c_2=c_0=0$ around each $e^{i\beta_l}$, it is easy to use maximum principle to prove $\tilde w_k=o(1)$ in $B_3$ using $\tilde w_k=o(1)$ on $\partial B_3$ and the Green's representation of $\tilde w_k$. The smallness of $\tilde w_k$ immediately implies $\tilde w_k=o(1)$ in $B_R$ for any fixed $R>>1$. Outside $B_R$, a crude estimate of $v_{0,k}k$ is
  $$v_{0,k}(y)\le -\mu_k-4(N+1)\log |y|+C, \quad 3<|y|<\tau \delta_k^{-1}. $$
  Using this and the Green's representation of $w_k$ we can first observe that the oscillation on each $\partial B_r$ is $o(1)$ ($R<r<\delta_k^{-1}/2$) and then by the Green's representation of $\tilde w_k$ and fast decay rate of $e^{V_k}$ we obtain $\tilde w_k=o(1)$ in $\overline{B(0, \delta_k^{-1})}$. A contradiction to $\max |\tilde w_k|=1$.

 There are $N+1$ local maximums. Correspondingly there are $N+1$ global solutions $V_{l,k}$ that
approximate $v_k$ accurately near $Q_l^k$ for $l=0,...,N$.  For $V_{l,k}$ the expression is
$$V_{l,k}=\log \frac{e^{\mu_l^k}}{(1+\frac{e^{\mu_l^k}h_{0,k}(Q_l^k)}{D}|y^{N+1}-(e_1+p_l^k)|^2)^2},\quad l=0,...,N, $$
where $p_l^k=E$ and
\begin{equation}\label{def-D}
D=8(N+1)^2.
\end{equation}
The equation that $V_{l,k}$ satisfies is 
$$\Delta V_{l,k}+|y|^{2N}h_{0,k}(Q_l^k)e^{V_{l,k}}=0,\quad \mbox{in}\quad \mathbb R^2. $$
Note that there is no need to define a $\phi_{l,k}$ and $v_{l,k}=v_k-\phi_{l,k}$ because the difference is insignificant.

Since $h_{0,k}(y)=e^{\phi_1}$ for some harmonic function $\phi_1(y)=O(\delta_k^*)$ for $|y|\sim 1$, we have 

\begin{equation}\label{extra-error}
h_{0,k}(Q_l^k)-h_{0,k}(Q_s^k)=O(\delta_k^*)=o(1)M_k. 
\end{equation}

Since $v_{0,k}$ and $V_{l,k}$ have the same common local maximum at $Q_l^k$, it is easy to see that
\begin{equation}\label{ql-exp}
Q_l^k=e^{i\beta_l}+\frac{p_l^ke^{i\beta_l}}{N+1}+O(|p_l^k|^2),\quad \beta_l=\frac{2l\pi}{N+2}.
\end{equation}
Let $M_{l,k}$ be the maximum of $|v_{0,k}-V_{l,k}|$ and we claim that all these $M_{l,k}$ are comparable:
\begin{equation}\label{M-comp}
M_{l,k}\sim M_{s,k},\quad \forall s\neq l.
\end{equation}
The proof of (\ref{M-comp}) is as follows: We use $L_{s,l}$ to denote the limit of $(v_{0,k}-V_{l,k})/M_{l,k}$ around $Q_s^k$:
$$\frac{(v_{0,k}-V_{l,k})(Q_s^k+\epsilon_kz)}{M_{l,k}}=L_{s,l}+o(1),\quad |z|\le \tau \epsilon_k^{-1} $$
where
$$ L_{s,l}=c_{1,s,l}\frac{z_1}{1+\frac 18 |z|^2}+c_{2,s,l}\frac{z_2}{1+\frac 18 |z|^2},\quad \mbox{and}\quad L_{l,l}=0, \quad s=0,...,N. $$
If all $c_{1,s,l}$ and $c_{2,s,l}$ are zero for a fixed $l$, we can obtain a contradiction just like the beginning of step two. So at least one of them is not zero.
For each $s\neq l$, by Lemma \ref{t-w-1-better} we have
\begin{equation}\label{Q-bad}
v_{0,k}(Q_s^k+\epsilon_kz)-V_{s,k}(Q_s^k+\epsilon_kz)=o(\epsilon_k)(1+|z|) M_{s,k},\quad |z|<\tau \epsilon_k^{-1}.
\end{equation}
Let $M_k=\max_{i}M_{i,k}$ ($i=0,...,N$) and we suppose $M_k=M_{l,k}$. Then to determine $L_{s,l}$ we see that
\begin{align*}
    &\frac{v_{0,k}(Q_s^k+\epsilon_k z)-V_{l,k}(Q_s^k+\epsilon_kz)}{M_k}\\
    =&o(\epsilon_k)(1+|z|)+\frac{V_{s,k}(Q_s^k+\epsilon_kz)-V_{l,k}(Q_s^k+\epsilon_kz)}{M_k}. 
\end{align*}
This expression says that $L_{s,l}$ is mainly determined by the difference of two global solutions $V_{s,k}$ and $V_{l,k}$. In order to obtain a contradiction to our assumption we will put the difference in several terms. The main idea in this part of the reasoning is that ``first order terms" tell us what the kernel functions should be, then the ``second order terms" tell us where the pathology is. 

We write $V_{s,k}(y)-V_{l,k}(y)$ as
$$V_{s,k}(y)-V_{l,k}(y)=\mu_s^k-\mu_l^k+2A-A^2+O(|A|^3) $$
where
$$A(y)=\frac{\frac{e^{\mu_l^k}}{D}|y^{N+1}-e_1-p_l^k|^2-\frac{e^{\mu_s^k}}{D}|y^{N+1}-e_1-p_s^k|^2}{1+\frac{e^{\mu_s^k}}{D}|y^{N+1}-e_1-p_s^k|^2}.$$
Here for convenience we abuse the notation $\epsilon_k$ by assuming $\epsilon_k=e^{-\mu_s^k/2}$. Note that $\epsilon_k=e^{-\mu_t^k/2}$ for some $t$, but it does not matter which $t$ it is. The difference between $h_{l,k}(Q_l^k)$ and $h_{l,k}(Q_s^k)$ is in (\ref{extra-error}). From $A$ we claim that

\begin{align}\label{late-1}
&V_{s,k}(Q_s^k+\epsilon_kz)-V_{l,k}(Q_s^k+\epsilon_kz)\\
=&\phi_1+\phi_2+\phi_3+\phi_4+\mathfrak{R},\nonumber
\end{align}
where
\begin{align}\label{4-express}
&\phi_1=(\mu_s^k-\mu_l^k)(1-\frac{h_s}{8}|z+\frac{N}2\epsilon_k z^2e^{-i\beta_s}|^2)/B, \\
&\phi_2=\frac{h_s}{2B}Re((z+\frac{N}2\epsilon_ke^{-i\beta_s}z^2))(\frac{\bar p_s^k-\bar p_l^k}{\epsilon_k}e^{-i\beta_s}))
\nonumber \\
&\phi_3=\frac{|p_s^k-p_l^k|^2}{4(N+1)^2\epsilon_k^2 B}\bigg (1-
\frac{|z|^2(1+2\cos (2\theta-2\theta_{st}-2\beta_s)}{8B} \bigg ),\nonumber\\
&\phi_4=\frac{h_s}4\frac{h_l-h_s}{h_l}|z+O(\epsilon_k)|z|^2|^2/B.\nonumber\\
&B=1+\frac {h_s}8|z+\frac N2z^2e^{-i\beta_s}\epsilon_k|^2, \nonumber
\end{align}
where $h_s=h_{0,k}(Q_s^k)$, $h_l=h_{0,k}(Q_l^k)$, $\mathfrak{R}_k$ is the collections of other insignificant terms.  Here we briefly explain the roles of each term. $\phi_1$ corresponds to the radial solution in the kernel of the linearized operator of the global equation. In other words, $\phi_1^k/M_k$ should tend to zero because in step one we have proved $c_0=0$. $\phi_2^k/M_k$ is the combination of the two other functions in the kernel. $\phi_3+\phi_4$ is the second order term which will play a leading role later.  The derivation of (\ref{late-1}) is put in the appendix.
 
Here $\phi_1$, $\phi_2$ correspond to solutions to the linearized operator. Here we note that if we set $\epsilon_{l,k}=e^{-\mu_l^k/2}$, there is no essential difference between $\epsilon_{l,k}$ and $\epsilon_k=e^{-\frac 12\mu_{1,k}}$ because $\epsilon_{l,k}=\epsilon_k(1+o(1))$. If $|\mu_{s,k}-\mu_{l,k}|/M_k\ge C$ there is no way to obtain a limit in the form of $L_{s,l}$ mentioned before. Thus we must have $|\mu_{s,k}-\mu_{l,k}|/M_k\to 0$. After simplification (see $\phi_2$ of (\ref{late-1})) we have
\begin{align}\label{c-12}
c_{1,s,l}=\lim_{k\to \infty}\frac{|p_s^k-p_l^k|}{2(N+1)M_k\epsilon_k}\cos(\beta_s+\theta_{sl}),\\
c_{2,s,l}=\lim_{k\to \infty}
\frac{|p_s^k-p_l^k|}{2(N+1)\epsilon_k M_k}sin(\beta_s+\theta_{sl})\nonumber
\end{align}
We omit $k$ for convenience.
It is also important to observe that even if $M_k=o(\epsilon_k)$ we still have $M_k\sim \max_{s}|p_s^k-p_l^k|/\epsilon_k$. Since each $|p_l^k|=E$, an upper bound for $M_k$ is 
\begin{equation}\label{small-M}
M_k\le C\mu_k\epsilon_k.
\end{equation}

Equation (\ref{c-12}) gives us a key observation: $|c_{1,s,l}|+|c_{2,s,l}|\sim |p_s^k-p_l^k|/(\epsilon_k M_k)$. So whenever $|c_{1,s,l}|+|c_{2,s,l}|\neq 0$ we have
$\frac{|p_s^k-p_l^k|}{\epsilon_k}\sim M_k$. In other words for each $l$, $M_{l,k}\sim \max_{t\neq l}\frac{|p_t^k-p_l^k|}{\epsilon_k}$.  Hence for any $t$, if $\frac{|p_t^k-p_l^k|}{\epsilon_k}\sim M_k$, let $M_{t,k}$ be the maximum of $|v_k-V_{t,k}|$, we have $M_{t,k}\sim M_k$. If all $\frac{|p_t^k-p_l^k|}{\epsilon_k}\sim M_k$ (\ref{M-comp}) is proved. So we prove that even if some $p_t^k$ is very close to $p_l^k$, $M_t^k$ is still comparable to $M_k$. The reason is there exists $q$ such that $\frac{|p_l^k-p_q^k|}{\epsilon_k}\sim M_k$, if $\frac{|p_t^k-p_l^k|}{\epsilon_k}=o(1)M_k$,
$$|p_t^k-p_q^k|\ge |p_l^k-p_q^k|-|p_t^k-p_l^k|\ge \frac 12 |p_l^k-p_q^k|.$$
Thus $\frac{|p_t^k-p_q^k|}{\epsilon_k}\sim M_k$ and $M_t^k\sim M_k$. (\ref{M-comp}) is established.  From now on for convenience we shall just use $M_k$, which has an upper bound in (\ref{small-M}).

Set $w_{l,k}=(v_{0,k}-V_{l,k})$, then we have $w_{l,k}(Q_l^k)=|\nabla w_{l,k}(Q_l^k)|=0$. Correspondingly we set
$$\tilde w_{l,k}=w_{l,k}/M_k. $$
The equation of $w_{l,k}$ can be written as 
\begin{equation}\label{wlk-bs}
   \Delta \tilde w_{l,k}+|y|^{2N}h_{0,k}(Q_l^k)e^{\xi_l^k}\tilde w_{l,k}=o(\sigma_k)(y-Q_l^k)|y|^{2N}e^{V_{l,k}}
   +o(\sigma_k)|y-Q_l^k|^2e^{V_{l,k}},
   \end{equation}
where $\xi_l^k$ comes from the Mean Value Theorem and satisfies
\begin{equation}\label{around-ls}
e^{\xi_l^k}=e^{V_{l,k}}(1+\frac 12w_{l,k}+O(w_{l,k}^2))=e^{V_{l,k}}(1+\frac 12M_k\tilde w_{l,k}+O(M_k^2)). 
\end{equation}
The function $\tilde w_{l,k}$ satisfies
\begin{equation}\label{aroud-s-1}
\lim_{k\to \infty}\tilde w_{l,k}(Q_s^k+\epsilon_k z)=\frac{c_{1,s,l}z_1+c_{2,s,l}z_2}{1+\frac 18 |z|^2}
\end{equation}
and around each $Q_s^k$ (\ref{Q-bad}) holds with $M_{s,k}$ replaced by $M_k$.

Now for $|y|\sim 1$, we use $\tilde w_{l,k}(Q_l^k)=0$ to write $\tilde w_{l,k}(y)$ as
\begin{align}\label{eva-wlk}
    \tilde w_{l,k}(y)&=\int_{\Omega_k}H_{y,l}(\eta)|\eta |^{2N}h_{0,k}(\eta)e^{\xi_l}\tilde w_{l,k}(\eta)d\eta \\
&+o(1)\delta_k^{N+1}e^{-\mu_k}. \nonumber
\end{align}
where the last term is based on (\ref{bry-small}) and
$$H_{y,l}(\eta):=G_k(y,\eta)-G_k(Q_l^k,\eta).$$
If we only consider $|y|\sim 1$, we use
\begin{equation}\label{crude-H}
H_{y,l}(\eta)=-\frac 1{2\pi}\log \frac{|y-\eta |}{|Q_l^k-\eta |}+O(\delta_k^2). 
\end{equation}
To evaluate the right hand side of (\ref{eva-wlk}), we only need to concentrate on $B(Q_s^k,\tau)$ for $s\neq l$ because the integration around $Q_l^k$ and outside bubbling disks only contribute $o(\epsilon_k)$. 
Around each $Q_s^k$ the $e^{\xi_l}$ can be replaced by $e^{V_{s,k}}$ with controllable error because from Lemma \ref{t-w-1-better}, the difference between them only leads to $o(\epsilon_k)$ in error. Also, $h_{0,k}(\eta)$ can be replaced by $h_{0,k}(Q_s^k)$ because
$$h_{0,k}(\eta)=h_{0,k}(Q_s^k)+\nabla h_{0,k}(Q_s^k)(\eta-Q_s^k)+O(\delta_k^*)(\eta-Q_s^k)^2. $$
Since $\nabla h_{0,k}(Q_s^k)=O(\delta_k^*)$, we see that the difference between $h_{0,k}(\eta)$ and $h_{0,k}(Q_s^k)$ only leads to $o(\epsilon_k)$ in its corresponding integration.

From the decomposition in (\ref{late-1}) we can now estimate the integral of $\tilde w_{l,k}$ more precisely. Clearly we only need to evaluate integrals around each $Q_s^k$. For this we have

\begin{equation}\label{late-2}
\int_{B(Q_s^k,\tau)}\tilde w_{l,k}(\eta)h_{0,k}(Q_l)|\eta |^{2N}e^{V_{s,k}}d\eta+o(\epsilon_k) 
=D_{s,l}^k+o(\epsilon_k).
\end{equation}
where
$$D_{s,l}^k=\frac{\pi}{2(N+1)^2}\frac{|p_s^k-p_l^k|^2}{\epsilon_k^2 M_{k}^2}M_k+2\pi\frac{h_l-h_s}{M_k}$$
and for convenience we skip $k$ in $Q_l$, and we shall write the two components of $Q_l$ as $Q_l=(Q_l^1,Q_l^2)$. 
The derivation of (\ref{late-2}) can be found in Appendix B.
Then for $|y|\sim 1$
\begin{align*}
    &\tilde w_{l,k}(y)=-\sum_{s\neq l}(H_{y,l}(Q_s)+O(\delta_k^2))D_{s,l}^k \\
    &-\sum_{s\neq l}\int_{B(Q_s,\tau)}\bigg ((\partial_1H_{y,l}(Q_s)\eta_1+\partial_2H_{y,l}(Q_s)\eta_2)h_{0,k}(Q_l)|\eta |^{2N}e^{\xi_l}\tilde w_{l,k}(\eta)\bigg )d\eta+o(\epsilon_k).
\end{align*} 

After evaluation we have
\begin{align*}
    \tilde w_{l,k}(y)&=-\frac{1}{2\pi}\sum_{s\neq l}(\log\frac{|y-Q_s|}{|Q_l-Q_s|}+O(\delta_k^2))D_{s,l}^k\\
   & -\sum_{s\neq l}8\Big(\frac{y_1-Q_s^1}{|y-Q_s|^2}-\frac{Q_l^1-Q_s^1}{|Q_l-Q_s|^2}\Big)c_{1,s,l}\epsilon_k\\
   &-8\Big(\frac{y_2-Q_s^2}{|y-Q_s|^2}-\frac{Q_l^2-Q_s^2}{|Q_s-Q_l|^2}\Big)c_{2,s,l}\epsilon_k+o(\epsilon_k).
\end{align*}
where we used
$$\int_{\mathbb R^2}\frac{z_1^2}{(1+\frac 18|z|^2)^3}dz=\int_{\mathbb R^2}\frac{z_2^2}{(1+\frac 18|z|^2)^3}dz=16\pi. $$
Recall that 
$$c_{1,s,l}=\frac{|p_s-p_l|}{2(N+1)M_k\epsilon_k}\cos (\beta_s+\theta_{sl})$$
$$c_{2,s,l}=\frac{|p_s-p_l|}{2(N+1)M_k\epsilon_k}\sin (\beta_s+\theta_{sl})$$

For $|y|\sim 1$ but away from the $N+1$ bubbling disks, we have, for $l\neq s$, 
$$v_{0,k}(y)=V_{l,k}(y)+M_k\tilde w_{l,k}(y) $$
and
$$v_{0,k}(y)=V_{s,k}(y)+M_k\tilde w_{s,k}(y). $$
Thus for $s\neq l$ we have
\begin{equation}\label{compare-10}
\frac{V_{s,k}(y)-V_{l,k}(y)}{M_k}=\tilde w_{l,k}(y)-\tilde w_{s,k}(y). \end{equation}

For $|y|\sim 1$ away from bubbling disks, we have
\begin{align*}&V_{s,k}(y)-V_{l,k}(y)\\
=&\mu_l^k-\mu_s^k+2\log \frac{h_s}{h_l}\\
&+2\log (1+\frac{\frac{D}{h_l}e^{-\mu_l^k}-\frac{D}{h_s}e^{-\mu_s^k}+ |y^{N+1}-1-p_l|^2-|y^{N+1}-1-p_s|^2}{\frac{D}{h_s}e^{-\mu_s^k}|y^{N+1}-e_1-p_s|}
\end{align*}
To evaluate the above, we use
\begin{align*} &|y^{N+1}-e_1-p_l|^2\\
=&|y^{N+1}-e_1-p_s|^2+2Re((y^{N+1}-e_1-p_s)(\bar p_s-\bar p_l))+|p_l-p_s|^2. 
\end{align*}
and set
$$\sigma_k:=(\mu_l^k-\mu_s^k+2\log \frac{h_s}{h_l})/M_k.$$
Then we see that $\sigma_k=o(1)$ and $\frac{D}{h_l}e^{-\mu_l^k}-\frac{D}{h_s}e^{-\mu_s^k}=o(M_k)\epsilon_k^2$. 
Thus
\begin{align*}
&\frac{V_{s,k}(y)-V_{l,k}(y)}{M_k}\\
=&4Re(\frac{y^{N+1}-1}{|y^{N+1}-1|^2}\frac{\bar p_l-\bar p_s}{M_k\epsilon_k}) \epsilon_k+\sigma_k+o(\epsilon_k). 
\end{align*}
On the other hand, for $y\in B_5\setminus (\cup_{t=1}^NB(Q_t^k,\tau_1))$, 
\begin{align*}
 &\tilde w_{l,k}(y)-\tilde w_{s,k}(y)\\
 =&-\frac{1}{2\pi}\sum_{m,m\neq l}\log \frac{|y-Q_m|}{|Q_l-Q_m|}D_{ml}^k\\
 &+8\epsilon_k\sum_{m,m\neq l}\bigg ((\frac{y_1-Q_m^1}{|y-Q_m|^2}-\frac{Q_l^1-Q_m^1}{|Q_l-Q_m|^2})\frac{|p_m-p_l|}{2(N+1)M_k\epsilon_k}\cos(\beta_m+\theta_{ml})\\
 &+(\frac{y_2-Q_m^2}{|y-Q_m|^2}-\frac{Q_l^2-Q_m^2}{|Q_l-Q_m|^2})\frac{|p_m-p_l|}{2(N+1)M_k\epsilon_k}\sin (\beta_m+\theta_{ml})\bigg )\\
 &+\frac{1}{2\pi}\sum_{m,m\neq s}\log \frac{|y-Q_m|}{|Q_l-Q_m|}D_{ms}^k\\
 &-8\epsilon_k\sum_{m,m\neq s}\bigg ((\frac{y_1-Q_m^1}{|y-Q_m|^2}-\frac{Q_s^1-Q_m^1}{|Q_s-Q_m|^2})\frac{|p_m-p_s|}{2(N+1)M_k\epsilon_k}\cos(\beta_m+\theta_{ms})\\
 &+(\frac{y_2-Q_m^2}{|y-Q_m|^2}-\frac{Q_s^2-Q_m^2}{|Q_s-Q_m|^2})\frac{|p_m-p_s|}{2(N+1)M_k\epsilon_k}\sin (\beta_m+\theta_{ms})\bigg )
\end{align*}
for all $l\neq s$. If we fix a set of $l,s$ that corresponds to the largest $|D_{ms}^k|$ and we consider $y$ close to $Q_s^k$. If we use $y=e^{i\beta_s}+z$ by abusing the notation $z$, then we have
\[y^{N+1}=(e^{i\beta_s}(1+ze^{-i\beta_s}))^{N+1}=1+(N+1)ze^{-i\beta_s}+O(|z|^2).\]
Therefore
\begin{align*}
&4Re(\frac{y^{N+1}-1}{|y^{N+1}-1|^2}\frac{\bar p_s-\bar p_l}{M_k\epsilon_k}\\
=&\frac{4|p_s-p_l|}{(N+1)|z|^2M_k\epsilon_k}\bigg (z_1\cos(\beta_s+\beta_{sl})+z_2\sin(\beta_s+\beta_{sl})+O(|z|^2)\bigg ).
\end{align*}
In the expression of $\tilde w_{l,k}(y)-\tilde w_{s,k}(y)$, we identify the leading term, which is 
\begin{align*}
8\epsilon_k\bigg ((\frac{y_1-Q_s^1}{|y-Q_s|^2}-\frac{Q_l^1-Q_s^1}{|Q_l-Q_s|^2})\frac{|p_s-p_l|}{2(N+1)M_k\epsilon_k}\cos(\beta_s+\theta_{sl})\\
+(\frac{y_2-Q_s^2}{|y-Q_s|^2}-\frac{Q_l^2-Q_s^2}{|Q_l-Q_s|^2})\frac{|p_s-p_l|}{2(N+1)M_k\epsilon_k}\sin(\beta_s+\theta_{sl})\bigg )\\
-\frac{1}{2\pi}\log \frac{|y-Q_s|}{|Q_l-Q_s|}D_{sl}^k.
\end{align*}
If we use $y=e^{i\beta_s}+z$ for $|z|$ small and replace $Q_s$ by $e^{i\beta_s}$ because their difference is $o(\epsilon_k)$. Then the expression above has this leading term:
\[\frac{4|p_s-p_l|}{(N+1)|z|^2M_k\epsilon_k}\bigg (z_1\cos(\beta_s+\beta_{sl})+z_2\sin(\beta_s+\beta_{sl})\bigg )-\frac 1{2\pi}\log \frac{|z|}{|e^{i\beta_l}-e^{i\beta_s}|}D_{sl}^k.\]
This we obtain $D_{sl}^k/\epsilon_k=o(1)$. Therefore 
\begin{equation}\label{dsl0}
D_{sl}^k=o(\epsilon_k),\quad \forall s\neq l.
\end{equation}
With this updated information we write $\tilde w_{l,k}(y)-\tilde w_{s,k}(y)$ as
\begin{align*}
 &\tilde w_{l,k}(y)-\tilde w_{s,k}(y)\\
 &=8\epsilon_k\sum_{m,m\neq l}\bigg ((\frac{y_1-\cos \beta_m}{|y-e^{i\beta_m}|^2}-\frac{\cos \beta_l-\cos \beta_m}{|e^{i\beta_l}-e^{i\beta_m}|^2})\frac{|p_m-p_l|}{2(N+1)M_k\epsilon_k}\cos(\beta_m+\theta_{ml})\\
 &+(\frac{y_2-\sin \beta_m}{|y-e^{i\beta_m}|^2}-\frac{\sin\beta_l-\sin \beta_m}{|e^{i\beta_l}-e^{i\beta_m}|^2})\frac{|p_m-p_l|}{2(N+1)M_k\epsilon_k}\sin (\beta_m+\theta_{ml})\bigg )\\
 &-8\epsilon_k\sum_{m,m\neq s}\bigg ((\frac{y_1-\cos \beta_m}{|y-e^{i\beta_m}|^2}-\frac{\cos \beta_s-\cos \beta_m}{|e^{i\beta_s}-e^{i\beta_m}|^2})\frac{|p_m-p_s|}{2(N+1)M_k\epsilon_k}\cos(\beta_m+\theta_{ms})\\
 &+(\frac{y_2-\sin \beta_m}{|y-e^{i\beta_m}|^2}-\frac{\sin \beta_s-\sin \beta_m}{|e^{i\beta_s}-e^{i\beta_m}|^2})\frac{|p_m-p_s|}{2(N+1)M_k\epsilon_k}\sin (\beta_m+\theta_{ms})\bigg )+o(\epsilon_k)
\end{align*}
for $|y|\in B_5\setminus (\cup_{l=1}^N B(Q_l^k,\tau_1))$. 

\medskip

Now in particular we take $l=0$ and we use the following notations: $\tilde w_k$, $V_k$, $c_{1,s}$, $c_{2,s}$, $\theta_s$, instead of $\tilde w_{0}^k$, $v_{0,k}$, $c_{1,s,0}$, $c_{2,s,0}$, $\theta_{s,0}$. 

The expression of $\tilde w_k$ gives
\begin{align*}
&\nabla \tilde w_k(y)\\
=&\int_{\Omega_k} \nabla_y G(y,\eta)\bigg (\mathfrak{h}_k(\delta_k Q_0)|\eta |^{2N}e^{\xi_k}\tilde w_k(\eta)
+\sigma_k\nabla\mathfrak{h}_k(\delta_k Q_0)(\eta-Q_0)|\eta |^{2N}e^{V_k(\eta)}\\
&+\frac{\delta_k^2}{M_k}\sum_{|\alpha |=2}\frac{\partial^{\alpha}\mathfrak{h}_k(\delta_k Q_0)(\eta -Q_0)^{\alpha}}{\alpha !}|\eta |^2{2N}e^{V_k(\eta)}\bigg )d\eta
+o(\epsilon_k), 
\end{align*}
for $y\in B_5\setminus (\cup_{s=1}^N B(Q_s^k,\tau_1))$.
Now we take $y=Q_0$, we have 
\begin{align*}
&0=\nabla \tilde w_k(Q_0)\\
=&\int_{\Omega_k} (-\frac{1}{2\pi})\frac{Q_0-\eta}{|Q_0-\eta |^2}\bigg (\mathfrak{h}_k(\delta_k Q_0)|\eta |^{2N}e^{\xi_k}\tilde w_k(\eta)
+\sigma_k\nabla\mathfrak{h}_k(\delta_k Q_0)(\eta-Q_0)|\eta |^{2N}e^{V_k(\eta)}\\
&+\frac{\delta_k^2}{M_k}\sum_{|\alpha |=2}\frac{\partial^{\alpha}\mathfrak{h}_k(\delta_k Q_0)(\eta -Q_0)^{\alpha}}{\alpha !}|\eta |^2{2N}e^{V_k(\eta)}\bigg )d\eta
+o(\epsilon_k), 
\end{align*}
Obviously we will integrate each of the two components in $B(Q_s^k,\tau_1)$ for $\tau_1>0$ small. Then we observe (\ref{dsl0}) 
that
\begin{align*}\int_{B(Q_l,\tau_1)}\bigg (\mathfrak{h}_k(\delta_kQ_0)|\eta |^{2N}e^{\xi_k}\tilde w_k(\eta)
+\sigma_k\nabla\mathfrak{h}_k(\delta_kQ_0)(\eta-Q_0)|\eta |^{2N}e^{V_k(\eta)}\\
+\frac{\delta_k^2}{M_k}\sum_{|\alpha |=2}\frac{\partial^{\alpha}\mathfrak{h}_k(\delta_k Q_0)(\eta -Q_0)^{\alpha}}{\alpha !}|\eta |^2{2N}e^{V_k(\eta)}\bigg )d\eta=o(\epsilon_k).\end{align*}
Based on this we use the following format: If $f$ is a smooth function,
\begin{align*}
&\int_{B(Q_s,\tau_1)} f(\eta) \bigg (\mathfrak{h}_k(\delta_k Q_0)|\eta |^{2N}e^{\xi_k}\tilde w_k(\eta)
+\sigma_k\nabla\mathfrak{h}_k(\delta_k Q_0)(\eta-Q_0)|\eta |^{2N}e^{V_k(\eta)}\\
&\qquad +\frac{\delta_k^2}{M_k}\sum_{|\alpha |=2}\frac{\partial^{\alpha}\mathfrak{h}_k(\delta_k Q_0)(\eta -Q_0)^{\alpha}}{\alpha !}|\eta |^2{2N}e^{V_k(\eta)}\bigg )d\eta\\
&=\partial_1f(e^{i\beta_s})c_{1s}\cdot 16\pi \epsilon_k+\partial_1 f(e^{i\beta_s})c_{2s}\cdot 16\pi \epsilon_k +o(\epsilon_k). 
\end{align*}
Then we replace $f(\eta_1,\eta_2)$ by \[f_1(\eta_1,\eta_2)=(-\frac{1}{2\pi})\frac{1-\eta_1}{(1-\eta_1)^2+\eta_2^2}\] and 
\[f_2(\eta_1,\eta_2)=\frac{1}{2\pi}\frac{\eta_2}{(1-\eta_1)^2+\eta_2^2}.\]
Then we have, from the expressions of $c_{1,s}$, $c_{2,s}$ in (\ref{c-12}), that 
\begin{align*}
0&=\partial_1\tilde w_k(Q_0)\\
&=16\pi\epsilon_k\sum_{s=1}^N\bigg (\frac{\cos(\beta_s+\theta_s)\cos \beta_s+\sin \beta_s\sin (\beta_s+\theta_s)}{4\pi(1-\cos \beta_s)}\frac{|p_s-p_0|}{2(N+1)M_k\epsilon_k}\bigg )+o(\epsilon_k)\\
&=4\epsilon_k\sum_{s=1}^N\frac{\cos \theta_s}{1-\cos \beta_s}\frac{|p_s-p_0|}{2(N+1)M_k\epsilon_k} +o(\epsilon_k). 
\end{align*}
Similarly
\begin{align*}
0&=\partial_2\tilde w_k(Q_0)\\
&=16\pi\epsilon_k\sum_{s=1}^N\bigg (\frac{\cos(\beta_s+\theta_s)\sin \beta_s-\cos \beta_s\sin (\beta_s+\theta_s)}{4\pi(1-\cos \beta_s)}\frac{|p_s-p_0|}{2(N+1)M_k\epsilon_k}\bigg )+o(\epsilon_k)\\
&=-4\epsilon_k\sum_{s=1}^N\frac{\sin \theta_s}{1-\cos \beta_s}\frac{|p_s-p_0|}{2(N+1)M_k\epsilon_k} +o(\epsilon_k). 
\end{align*}
If we use $a_s$ to denote
\[a_s:=\lim_{k\to \infty}\frac{|p_s^k-p_0^k|}{(1-cos\beta_s)M_k\epsilon_k},\quad s=1,...,N,\]
then we have
\begin{equation}\label{eq-ar-1}
\sum_{s=1}^Na_s\cos \theta_s=0
\end{equation}
\begin{equation}\label{eq-ar-2}
\sum_{s=1}^Na_s\sin \theta_s=0,
\end{equation}
where $a_s\ge 0$ for all $s$. Taking the sum of the squares of (\ref{eq-ar-1}) and (\ref{eq-ar-2}) we obtain 
\[\sum_{s=1}^Na_s^2+\sum_{s<t}2a_sa_t\cos(\theta_s-\theta_t)=0. \]
Since $\sum_sa_s^2\ge 2\sum_{s<t}a_sa_t$, we have
\[\sum_{s<t}2a_sa_t(1+\cos(\theta_s-\theta_t))\le 0.\]
Since each term on the left is obviously non-negative, we know each 
\[a_sa_t(1+\cos (\theta_s-\theta_t))=0,\quad \forall s<t. \] 
If there is only one $a_s>0$, it is easy to see that (\ref{eq-ar-1}) and (\ref{eq-ar-2}) cannot both hold. If there are three $a_t's>0$, it is also elementary to see this is not possible: say $a_1,a_2,a_3>0$, then they have to be equal. Then we see that we must have 
\[\theta_1-\theta_2=\pm \pi,\quad \theta_2-\theta_3=\pm \pi,\quad \theta_1-\theta_3=\pm \pi .\] 
Obviously these three equations cannot hold at the same time. So the only situation left is there are exactly two $a_t's$ positive. All other $a_ts$ are zero. This means there are exactly two $p_{s_1}^k$, $p_{s_2}^k$ such that 
\[\lim_{k\to \infty}\frac{p_{s_1}^k-p_0^k}{\epsilon_kM_k}=-\lim_{k\to \infty}\frac{p_{s_2}^k-p_0^k}{\epsilon_k M_k}\neq 0,\quad \lim_{k\to \infty}\frac{p_t^k-p_0^k}{\epsilon_kM_k}=0,\quad \forall t\neq s_1,s_2.\]

If we apply the same argument to $\tilde w_l^k$. Then from $\nabla \tilde w_l^k(Q_l^k)=0$ we would get exactly $p_{l_1}^k$ and $p_{l_2}^k$ different from $p_l^k$ and 
\[\lim_{k\to \infty} \frac{p_{l_1}^k-p_l^k}{\epsilon_k M_k}=-\lim_{k\to \infty}\frac{p_{l_2}^k-p_l^k}{\epsilon_k M_k}\neq 0,\quad \lim_{k\to \infty}\frac{p_t^k-p_l^k}{\epsilon_k M_k}=0, \forall t\neq l_1,l_2.\]
Then it is easy to see that this is only possible when we have $N=2$ because if $N=1$, we would have just one $a_s\neq 0$, which is not possible based on (\ref{eq-ar-1}) and (\ref{eq-ar-2}). If $N\ge 3$, we have to have $p_t^k$ that satisfies 
\[\lim_{k\to \infty}\frac{|p_t^k-p_0^k|}{\epsilon_k M_k}=0,\quad \mbox{and}\quad \lim_{k\to \infty} \frac{p_t^k-p_{s_1}^k}{\epsilon_kM_k}=0\]
which is not possible because $\lim_{k\to \infty}\frac{p_{s_1}^k-p_0^k}{\epsilon_kM_k}\neq 0$. 

Finally we rule out the case $N=2$. In this case we have  \[\lim_{k\to \infty}\frac{p_1^k-p_0^k}{\epsilon_kM_k}=-\lim_{k\to \infty}\frac{p_2^k-p_0^k}{\epsilon_kM_k}\neq 0.\] 
However from $\tilde w_1^k(Q_1^k)=0$ we have 
\[\lim_{k\to \infty}\frac{p_2^k-p_1^k}{\epsilon_kM_k}=-\lim_{k\to \infty}\frac{p_0^k-p_1^k}{\epsilon_kM_k}\neq 0,\]
Then from $\tilde w_2^k(Q_2^k)=0$ we have 
\[\lim_{k\to \infty}\frac{p_1^k-p_2^k}{\epsilon_kM_k}=-\lim_{k\to \infty}\frac{p_0^k-p_2^k}{\epsilon_kM_k}\neq 0,\]
Easy to see it is not possible for all these three equations to hold at the same time.  
Lemma \ref{small-other} is established. $\Box$

\medskip

Proposition \ref{key-w8-8} is an immediate consequence of Lemma \ref{small-other}.  $\Box$.

\section{Proof of the main theorems}

First we make a simple observation. If we let 
$$\hat w_k=w_{0,k}/\delta_k^*. $$
It is easy to see that $\hat w_k$ tends to a harmonic function away from the $N+1$ local maximums of $v_{0,k}$. Since $\hat w_k$ is bounded, the global harmonic function $\hat w_k$ tends to has to be a constant. If we focus on the neighborhood of $Q_0^k$, by standard Fourier analysis as before we see that this constant is zero. Thus over compact subset away from the $N+1$ local maximum points of $v_{0,k}$ we have $\hat w_k=o(1)$ and $\nabla \hat w_k=o(1)$, which means in this region $w_{0,k}=o(\delta_k^*)$ and $\nabla w_{0,k}=o(\delta_k^*)$. Then we use the smallness of $w_{0,k}$ on $\partial B(Q_0,2r)$ to prove 
\begin{equation}\label{w0ksmall}
|w_{0,k}(Q_0^k+\epsilon_kz)|\le C\delta_k^*\epsilon_k(1+|z|),\quad |z|\le 2r\epsilon_k^{-1}. 
\end{equation}

In the next step we consider the difference between two Pohozaev identities over $\Omega_{0,k}:=B(Q_0,r)$
  for small $r>0$. For $v_k$ we have
\begin{align}\label{pi-vk}
\int_{\Omega_{0,k}}\partial_{\xi}(|y|^{2N}h_{0,k}(y))e^{v_{0,k}}-\int_{\partial \Omega_{0,k}}e^{v_{0,k}}|y|^{2N}h_{0,k}(y)(\xi\cdot \nu)\\
=\int_{\partial \Omega_{0,k}}(\partial_{\nu}v_{0,k}\partial_{\xi}v_{0,k}-\frac 12|\nabla v_{0,k}|^2(\xi\cdot \nu))dS. \nonumber
\end{align}
where $\xi$ is an arbitrary unit vector. Correspondingly the Pohozaev identity for $V_{0,k}$ is

\begin{align}\label{pi-Vk}
\int_{\Omega_{0,k}}\partial_{\xi}(|y|^{2N}h_{0,k}(Q_0^k))e^{V_{0,k}}-\int_{\partial \Omega_{0,k}}e^{V_k}|y|^{2N}h_{0,k}(Q_0^k)(\xi\cdot \nu)\\
=\int_{\partial \Omega_{0,k}}(\partial_{\nu}V_k\partial_{\xi}V_{0,k}-\frac 12|\nabla V_{0,k}|^2(\xi\cdot \nu))dS. \nonumber
\end{align}
Using the smallness of $w_{0,k}$ on $\partial \Omega_{0,k}$ we have

\begin{align*}&\int_{\partial \Omega_{s,k}}(\partial_{\nu}v_{0,k}\partial_{\xi}v_{0,k}-\frac 12|\nabla v_{0,k}|^2(\xi\cdot \nu))dS\\
-&\int_{\partial \Omega_{s,k}}(\partial_{\nu}V_{0,k}\partial_{\xi}V_{0,k}-\frac 12|\nabla V_{0,k}|^2(\xi\cdot \nu))dS
=o(\delta_k^*).
\end{align*}
Next we observe that the difference on the second terms of the two Pohozaev identities is minor: 
\[\int_{\partial \Omega_{0,k}}e^{v_{0,k}}|y|^{2N}h_{0,k}(y)(\xi\cdot \nu)-\int_{\partial \Omega_{0,k}}h_{0,k}(Q_0^k)e^{V_{0,k}}|y|^{2N}(\xi\cdot \nu)=o(\delta_k^*). \]

Finally by (\ref{w0ksmall}) we have
\begin{align}\label{imp-1}
&\int_{\Omega_{0,k}}\partial_{\xi}(|y|^{2N}h_{0,k}(y))e^{v_{0,k}}\\
=&\int_{\Omega_{0,k}}\partial_{\xi}(|y|^{2N}(h_{0,k}(Q_s^k)+\nabla h_{0,k}(Q_s^k)(y-Q_s^k))e^{V_{0,k}}(1+w_{0,k})dy+o(\delta_k^*)\nonumber\\
=&\int_{\Omega_{0,k}}(\partial_{\xi}(|y|^{2N})h_{0,k}(Q_s^k)e^{V_{0,k}}+\partial_{\xi}h_{0,k}(Q_s^k)|y|^{2N}e^{V_{0,k}}
+o(\delta_k^*).\nonumber
\end{align}
Thus we have proved $\nabla h_{0,k}(Q_0)=o(\delta_k^*)$. In a similar way we also obtain $|\nabla h_{0,k}(Q_s)|=o(\delta_k^*)$ for all $s=0,...,N$. Thus we get a contradiction because at least one of them is comparable to $\delta_k^*$. 
Theorem \ref{main-thm-1} is established. $\Box$

\bigskip

\noindent{\bf Proof of Theorem \ref{main-thm-gen}}: First it is a well known fact that all the blowup points are in the interior of $\Omega$ (see \cite{wei-ma}). Let $p$ be a blowup point of $\{\mathfrak{u}_k\}$, if $p$ is a regular point or a non-quantized singular point, the blowup solutions are known to be simple \cite{bclt,zhangccm}. So we only consider the case that $q$ is a quantized singular source. Obviously we can use a fundamental solution to reduce the study to equation (\ref{main-0}). Because there is at least one other blowup point, we can see that around $p$, the harmonic function that eliminates the oscillation of $u_k$ around $p$ does not tend to $0$, because around the other blowup point, the harmonic function is tending to infinity.  Theorem \ref{main-thm-gen} is an immediate consequence of Theorem \ref{main-thm-1}. $\Box$

\bigskip

\noindent{\bf Proof of Theorem \ref{conjecture}:}
We consider the blowup solutions of $u_k+\log \lambda_k$. It is well known that if $\lambda_k$ does not tend zero, $u_k$ is uniformly bounded. So we use $\hat u_k$ to denote $u_k+\log \lambda_k$ and we assume that $\hat u_k$ is a sequence of blowup solutions. If there are at least two blowup points, the argument for Theorem \ref{main-thm-gen} can be applied here to say that each blowup point is simple. So we only consider the case that there is one blowup point. Without loss of generality we take this blowup point to be the origin. If the origin is a regular point or a singular point with non-quantized Dirac mass, the blowup sequence is well known to be simple. Thus we assume that at $0$ there is a quantized singular source $4\pi N\delta_0$. If the domain is a disk centered at the origin, the proof of Theorem \ref{main-thm-1} also applies in this case that the only blowup point is simple. So the only case is that $\Omega$ is not a disk centered around the origin. In this case we assume that $B_1\subset \Omega\subset B_R$ for some $R>1$. 
Now we recall the asymptotic expansion of 
$$V_k(y)=\mu_k-2\log (1+\frac{e^{\mu_k}}{D}|y^{N++1}-1-p_k|^2) $$
where $|p_k|=O(\mu_ke^{-\mu_k})$. Thenf or $r$ large, say, comparable to $\delta_k^{-1}$, 
\begin{align*}V_k(y)=-\mu_k+2\log D-4(N+1)\log r
+\frac{4}{r^{N+1}}\cos((N+1)\theta)\\
+O(\delta_k^{N+1}\mu_ke^{-\mu_k})+O(\delta_k^{2N+2}). \end{align*}
We define a harmonic function $\phi_{0,k}$ that eliminates all the oscillation of $V_k$ on $\partial \Omega_k$. Consider the Fourier expansion of $\phi_{0,k}$:
$$\sum_{n=1}^{\infty}(a_n^k\cos(n\theta)+b_n^k\sin(n\theta))r^n. $$
Note that  there is no projection on $1$ because $\phi_{0,k}(0)=0$. Since the value of the absolute value of $\phi_{0,k}$ on $\partial \Omega_k$ (which is between $|y|=\delta_k^{-1}$ and $|y|=R\delta_k^{-1}$ ) is comparable to $\delta_k^{N+1}$, we make two observations. First the restriction of $\phi_{0,k}$ on $\partial B(0,\delta_k^{-1})$ is between $-C\delta_k^{N+1}$ and $C\delta_k^{N+1}$. This is because of the asymptotic expansion of $V_k$ for $|y|\sim \delta_k^{-1}$. This provides an upper bound of $a_n^k$ and $b_n^k$:
$$|a_n^k|+|b_n^k|\le C\delta_k^{N+1+n}. $$
Second, we shall prove there is a point $p_k\in \partial B(0,\delta_k^{-1})$  such that 
\begin{equation}\label{trivial-low}
|\phi_{0,k}(p_k)|\ge C\delta_k^{N+1}, \quad |p_k|=\delta_k^{-1}. 
\end{equation}

 Without loss of generality we can assume that $\partial \Omega_k$ is tangent to $\partial B(0,\delta_k)$. 
It is more convenient to consider this problem before scaling: Let $\phi$ be a harmonic function on $\Omega$
with $\phi(0)=0$ and the value of $\phi$ on $\partial \Omega$ can be described as
$$\phi(x)=f(|x|)\cos((N+1)\theta), \quad x=|x|e^{i\theta}\in \partial \Omega, $$
where $f(r)\sim 1$ for $1<r<R$. 
Let $e^{i\theta_0}$ be the tangent point between $\partial B_1$ and $\partial \Omega$, if $\cos((N+1)\theta_0\neq 0$, we already have 
$\phi(e^{i\theta_0})\neq 0$, which certainly gives
$|\phi(e^{i\theta_0})|\ge C$ for some $C>0$ independent of $k$. Once this is done, after scaling and normalization with $\delta_k^{N+1}$ we have (\ref{trivial-low})
for $p_k=\delta_k^{-1}e^{i\theta_0}$ and 
$\phi_{0,k}(y)=\delta_k^{N+1}\phi(\delta_k y)$. 

If $\cos((N+1)\theta_0)=0$, we just consider a neighborhood of $e^{i\theta_0}$. If $\theta$ is just slightly different from $\theta_0$, the magnitude of $\phi$ on $|x|e^{i\theta}\in \partial \Omega$ is comparable to $1$ because $f(|x|)$ is comparable to $1$. The curvature of $\partial \Omega$ is bounded, it is easy to use the Green's representation formula of $\phi$ to find $x_1=e^{i\theta_1}\in \partial B_1$, close to $e^{i\theta_0}$ such that $|\phi(x_1)|\ge C$ for some $C>0$. Scaling back to $\partial \Omega_k$ we have $p_k=\delta_k^{-1}e^{i\theta_1}$ and (\ref{trivial-low}) holds.

Based on these two observations we set 
$$\delta_k^*=\sum_{n=1}^{N+1}(|a_n^k|+|b_n^k|), $$
we have 
$$C\delta_k^{2N+2}\le \delta_k^*\le C\delta_k^{N+2}. $$
By setting $h_{0,k}(y)=e^{\phi_{0,k}}$ as before, we have, using the same proof of Lemma \ref{key-lem-2}, that 
$$|\nabla h_{0,k}(Q_s^k)|\ge C\delta_k^*, $$
for at least one $Q_s^k$. Then the remaining part of the proof is similar to that of Theorem \ref{main-thm-gen}. Theorem \ref{conjecture} is established. $\Box$.

\section{Appendix A}  
In this section we derive (\ref{late-1}) and (\ref{4-express}).
We use simplified notations:
$$V_s=\mu_s-2\log (1+\frac{e^{\mu_s}}{D_s}|y^{N+1}-1-p_s|^2).$$
$$V_l=\mu_l-2\log (1+\frac{e^{\mu_l}}{D_l}|y^{N+1}-1-p_l|^2).$$
where $D_s=\frac{8(N+1)^2}{h_s}$, $h_s=h_{0,k}(Q_s^k)$. 
$$V_s-V_l=\mu_s-\mu_l+2\log (1+A)=\mu_s-\mu_l+2A-A^2+O(A^3), $$
where
$$A=\frac{\frac{e^{\mu_l}}{D_l}|y^{N+1}-1-p_l|^2-\frac{e^{\mu_s}}{D_s}|y^{N+1}-1-p_s|^2}{1+\frac{e^{\mu_s}}{D_s}|y^{N+1}-1-p_s|^2} $$

\begin{align}\label{tem-5}
    \frac{e^{\mu_l}}{D_l}&=\frac{e^{\mu_s}}{D_s}e^{\mu_l-\mu_s}\frac{D_s}{D_l}\\
    &=\frac{e^{\mu_s}}{D_s}(1+(\mu_l-\mu_s)+O(\mu_l-\mu_s)^2)(1+\frac{h_l-h_s}{h_s}+O(\delta_k^{2*}))\nonumber\\
   & =\frac{e^{\mu_s}}{D_s}(1+\mu_l-\mu_s+\frac{h_l-h_s}{h_s}+O((\mu_l-\mu_s)^2+\delta_k^{2*})).\nonumber
\end{align}
where we used $(D_s-D_l)/D_l=(h_l-h_s)/h_s=O(\delta_k^*)$ and
we use $\delta_k^{2*}$ to denote $(\delta_k^*)^2$. 

Using this we now write $A$ as
\begin{align*}
A&=\bigg (\frac{e^{\mu_s}}{D_s}(|y^{N+1}-1-p_l|^2-|y^{N+1}-1-p_s|^2)\\
&+\frac{e^{\mu_s}}{D_s}(\mu_l-\mu_s+\frac{h_l-h_s}{h_s}+E)|y^{N+1}-1-p_l|^2\bigg )/\bigg (1+\frac{e^{\mu_s}}{D_s}|y^{N+1}-1-p_s|^2 \bigg )
\end{align*}
where $E=O((\mu_l-\mu_s)^2+\delta_k^{2*})$.

The derivation of (\ref{late-1}) is as follows: Recall that the expression of $Q_s$ is in (\ref{ql-exp}):
$$Q_s=e^{i\beta_s}(1+\frac{p_s}{N+1})+O(|p_s|^2).$$
and $y=Q_s+\epsilon_k z$. Direct computation shows
$$y^{N+1}=1+p_s+(N+1)\epsilon_kze^{-i\beta_s}+\frac{N(N+1)}2\epsilon_k^2z^2e^{-2i\beta_s}+O(\epsilon_k^3\mu_k|z|).$$
Then
$$|y^{N+1}-1-p_s|^2=(N+1)^2\epsilon_k^2|z+\frac{N}2\epsilon_kz^2e^{-i\beta_s}|^2+O(\epsilon_k^4\mu_k)|z|^2 $$
$$|y^{N+1}-1-p_l|^2=(N+1)^2\epsilon_k^2|z+\frac{N}2\epsilon_kz^2e^{-i\beta_s}+\frac{(p_s-p_l)e^{i\beta_s}}{(N+1)\epsilon_k }|^2+O(\epsilon_k^4\mu_k|z|^2). $$

With these expressions one can check easily that the expression of $B$ in (\ref{late-1}) holds:
$$B=1+\frac{h_s}{8}|z+\frac{N}2\epsilon_kz^2e^{-i\beta_s}|^2+O(\epsilon_k^2\mu_k)|z|^2. $$
Then the expression of $A$ is
\begin{align*}
    A=\frac{h_s}8\bigg ( 2Re(z+\frac{N}2\epsilon_kz^2e^{-i\beta_s})\frac{\bar p_s-\bar p_l}{(N+1)\epsilon_k}e^{-i\beta_s}+\frac{|p_s-p_l|^2}{(N+1)^2\epsilon_k^2}\\
    +(\mu_l-\mu_s+\frac{h_l-h_s}{h_s}+E_{\mu})|z+\frac{N}2\epsilon_kz^2e^{-i\beta_s}+O(\epsilon_k^2\mu_kz)|^2\bigg )/B.
\end{align*}
where $E_{\mu}=O(|\mu_l-\mu_s|^2+\delta_k^{2*})$. Note that $D_s=o(M_k)$, so $E_{\mu}=o(M_k^2)$.

\begin{equation}\label{exp-a-squ}
A^2=
\frac{1}{16(N+1)^2}\bigg (Re(z\frac{\bar p_s-\bar p_l}{\epsilon_k}e^{-i\beta_s})\bigg )^2/B^2
+\mbox{other terms}.
\end{equation}
The numerator of $A^2$ has the following leading term:
$$\frac{1}{32(N+1)^2}\bigg (|z|^2(\frac{|p_s-p_l|}{\epsilon_k})^2\big (1+2\cos(2\theta-2\theta_{sl}-2\beta_s)\big ) \bigg )$$
where $z=|z|e^{i\theta}$, $p_s-p_l=|p_s-p_l|e^{i\theta_{sl}}$.
Using these expressions we can obtain (\ref{late-1}) by direct computation.

\section{Appendix B}
In this section we prove (\ref{late-2}) based on (\ref{4-express}).
The terms of $\phi_1$ and $\phi_2$ lead to $o(\epsilon_k)$
The integrations involving $\phi_3$ and $\phi_4$ provide the leading term. More detailed information is the following:
First for a global solution
$$V_{\mu,p}=\log \frac{e^{\mu}}{(1+\frac{e^{\mu}}{\lambda}|z^{N+1}-p|^2)^2}$$ of
$$\Delta V_{\mu,p}+\frac{8(N+1)^2}{\lambda}|z|^{2N}e^{V_{\mu,p}}=0,\quad \mbox{in }\quad \mathbb R^2, $$
by differentiation with respect to $\mu$ we have
$$\Delta(\partial_{\mu}V_{\mu,p})+\frac{8(N+1)^2}{\lambda}|z|^{2N}e^{V_{\mu,p}}\partial_{\mu}V_{\mu,p}=0,\quad \mbox{in}\quad \mathbb R^2. $$
By the expression of $V_{\mu,p}$ we see that 
$$\partial_r\bigg (\partial_{\mu}V_{\mu,p}\bigg )(x)=O(|x|^{-2N-3}).$$ 
Thus we have
\begin{equation}\label{inte-eq-1}
\int_{\mathbb R^2}\partial_{\mu}V_{\mu,p}|z|^{2N}e^{V_{\mu,p}}=\int_{\mathbb R^2}\frac{(1-\frac{e^{\mu}}{\lambda}|z^{N+1}-P|^2)|z|^{2N}}{(1+\frac{e^{\mu}}{\lambda}|z^{N+1}-P|^2)^3}dz=0.
\end{equation}

From $V_{\mu,p}$  we also have
$$\int_{\mathbb R^2}\partial_{P}V_{\mu,p}|y|^{2N}e^{V_{\mu,p}}=\int_{\mathbb R^2}\partial_{\bar P}V_{\mu,p}|y|^{2N}e^{V_{\mu,p}}=0, $$
which gives
\begin{equation}\label{inte-eq-2}
\int_{\mathbb R^2}\frac{\frac{e^{\mu}}{\lambda}(\bar z^{N+1}-\bar P)|z|^{2N}}{(1+\frac{e^{\mu}}{\lambda}|z^{N+1}-P|^2)^3}=\int_{\mathbb R^2}\frac{\frac{e^{\mu}}{\lambda}( z^{N+1}- P)|z|^{2N}}{(1+\frac{e^{\mu}}{\lambda}|z^{N+1}-P|^2)^3}=0.
\end{equation}

Based on (\ref{inte-eq-1}) and (\ref{inte-eq-2}), the expressions of $\phi_1$, $\phi_2$ and $B$ in (\ref{4-express}) lead to 

$$\int_{B(0,\tau\epsilon_k^{-1})}\frac{\phi_1}{M_k}B^{-2}=o(\epsilon_k),\quad
\int_{B(0,\tau\epsilon_k^{-1})}\frac{\phi_2}{M_k}B^{-2}=o(\epsilon_k).$$
The integrations involving $\phi_3$ and $\phi_4$ lead to the expression of $D_{s,l}^k$ in (\ref{late-2}). 
Thus (\ref{late-2}) holds.

\section{Appendix C: Location of local maximums of $v_k$}

In this section we use standard arguments to prove that the local maximum points of $v_k$ is evenly distributed around the origin. Recall that $v_k$ satisfies (\ref{e-f-vk}). Here we abuse the notation by assuming that $Q_1^k=e_1$ is one local maximum of $v_k$ and we use $\mu_k$ instead of $\bar \mu_k$ to denote $v_k(e_1)$ for convenience.

Now we consider $v_k$ around $Q_l^k$. Using the results in \cite{chenlin1,zhangcmp,gluck} we have, for $v_k$ in $B(Q_l^k,\epsilon )$, the following gradient estimate:

 \begin{equation}\label{gra-each-p}
2N\frac{\tilde Q_l^k}{|\tilde Q_l^k|^2}+\nabla \phi_l^k(\tilde Q_l^k)=O( \mu_k e^{- \mu_k}),
\end{equation}
where $\phi_l^k$ is the harmonic function that eliminates the oscillation of $v_k$ on $\partial B(Q_l^k,\epsilon)$ and $\tilde Q_l^k$ is the maximum of $v_k-\phi_l^k$ that satisfies
\begin{equation}\label{close-2}
\tilde Q_l^k-Q_l^k=O( e^{-\mu_k}).
\end{equation}
Using (\ref{close-2}) in (\ref{gra-each-p}) we have
\begin{equation}\label{pi-each-p}
2N\frac{Q_l^k}{| Q_l^k|^2}+\nabla \phi_l^k( Q_l^k)=O( \mu_k e^{- \mu_k}).
\end{equation}
The first estimate of $\nabla \phi_l^k(Q_l^k)$ is
\begin{lem}\label{phi-k-e} For $l=0,...,N$,
\begin{equation}\label{gra-each-2}
\nabla \phi_l^k(Q_l^k)
=-4\sum_{m=0,m\neq l}^N\frac{Q_l^k-Q_m^k}{|Q_l^k-Q_m^k|^2}
+E
\end{equation}
where
$$E=O(\mu_k e^{-\mu_k}). $$
\end{lem}

\noindent{\bf Proof of Lemma \ref{phi-k-e}:}

\medskip

From the expression of $v_k$ on $\Omega_k=B(0,\tau \delta_k^{-1})$ we have, for $y$ away from bubbling disks,
\begin{align}\label{pi-e-1}
v_k(y)&=v_k|_{\partial \Omega_k}+\int_{\Omega_k}G(y,\eta)|\eta |^{2N}e^{v_k(\eta)}d\eta\\
&=v_k|_{\partial \Omega_k}+\sum_{l=0}^{N}G(y,Q_l^k)\int_{B(Q_l,\epsilon)}|\eta |^{2N}e^{v_k}d\eta \nonumber\\
&+\sum_l\int_{B(Q_l,\epsilon)}(G(y,\eta)-G(y,Q_l^k))|\eta |^{2N}e^{v_k}d\eta+O(\mu_ke^{-\mu_k}). \nonumber
\end{align}

We have, by standard estimates
$$
v_k(y)
=v_k|_{\partial \Omega_k}-4\sum_{l=0}^{N}\log |y- Q_l^k|
+8\pi\sum_{l=0}^{N}H(y,Q_l^k)+O(\mu_ke^{-\mu_k}). $$

The harmonic function that kills the oscillation of $v_k$ around $Q_m^k$ is
\begin{align*}
\phi_m^k=-4\sum_{l=0,l\neq m}^N(\log |y-Q_l^k|-\log |Q_m^k-Q_l^k|)\\
+8\pi\sum_{l=0}^{N}(H(y,Q_l^k)-H(Q_m^k,Q_l^k))+O(\mu_ke^{-\mu_k}).
\end{align*}

The corresponding estimate for $\nabla \phi_m^k$ is
$$
\nabla \phi_m^k(Q_m^k)=-4\sum_{l=0,l\neq m}^N\frac{Q_m^k-Q_l^k}{|Q_m^k-Q_l^k|^2}
+8\pi\sum_{l=0}^{N}\nabla_1 H(Q_m^k,Q_l^k)+O(\mu_ke^{-\mu_k}).
$$
where $\nabla_1$ stands for the differentiation with respect to the first component. From the expression of $H$, we have
\begin{align}\label{grad-H}
\nabla_1H(Q_m^k,Q_l^k)&=\frac 1{2\pi}\frac{Q_m^k-\tau^2\delta_k^{-2}Q_l^k/|Q_l^k|^2}{|Q_m^k-\tau^2\delta_k^{-2}Q_l^k/|Q_l^k|^2|^2}\\
&=\frac 1{2\pi}\tau^{-2}\delta_k^2\frac{\tau^{-2}\delta_k^2Q_m^k-Q_l^k/|Q_l^k|^2}{|Q_l^k/|Q_l^k|^2-\tau^{-2}\delta_k^2Q_m^k|^2}\nonumber \\
&=-\frac{1}{2\pi}\tau^{-2}\delta_k^2e^{\frac{2\pi i l}{N+1}}+O(\sigma_k\delta_k^2). \nonumber
\end{align}
where $\sigma_k=\max_l|Q_l^k-e^{\frac{2\pi i l}{N+1}}|$. Later we shall obtain more specific estimate of $\sigma_k$.

Thus
\begin{align}\label{phimk-g}
&\nabla \phi_m^k(Q_m^k) \\
=&-4\sum_{l=0,l\neq m}^N\frac{Q_m^k-Q_l^k}{|Q_m^k-Q_l^k|^2}
-4\tau^{-2}\delta_k^2\sum_{l=0}^{N}e^{\frac{2\pi il}{N+1}}+O(\sigma_k\delta_k^2)+O(\mu_ke^{-\mu_k}) \nonumber \\
=&-4\sum_{l=0,l\neq m}^N\frac{Q_m^k-Q_l^k}{|Q_m^k-Q_l^k|^2}
+O(\sigma_k\delta_k^2)+O(\mu_ke^{-\mu_k}) \nonumber
\end{align}
where we have used $\sum_{l=0}^{N}e^{2\pi l i/(N+1)}=0$.
Since we don't have the estimate of $\sigma_k$ now we have
$$
\nabla \phi_m^k(Q_m^k)=-4\sum_{l=0,l\neq m}^N\frac{Q_m^k-Q_l^k}{|Q_m^k-Q_l^k|^2}+
O(\mu_ke^{-\mu_k})+O(\sigma_k\delta_k^2).
$$
Lemma \ref{phi-k-e} is established. $\Box$

\bigskip

Using $$E_1=O(\mu_ke^{-\mu_k})+O(\sigma_k\delta_k^2). $$
The Pohozaev identity around $Q_l^k$ now reads

$$-4\sum_{j=0,j\neq l}^N\frac{Q_l^k-Q_j^k}{|Q_l^k-Q_j^k|^2}+2N\frac{Q_l^k}{|Q_l^k|^2}=E_1. $$
We have, treating every term as a complex number,
$$N\frac{1}{\bar Q_l^k}=2\sum_{j=0,j\neq l}^N\frac{1}{\bar Q_l^k-\bar Q_j^k}+E_1, $$
where $\bar Q_l^k$ is the conjugate of $Q_l^k$.
Thus
\begin{equation}\label{pi-N}
N=2\sum_{j=0,j\neq l}^N\frac{Q_l^k}{Q_l^k-Q_j^k}+E_1.
\end{equation}
Let $\beta_l=2\pi l/(N+1)$, we write $Q_l^k=e^{i\beta_l }+p_l^k$ for $p_l^k\to 0$. Then we write the first term on the right hand side of (\ref{pi-N}) as
\begin{align*}
&\frac{Q_l^k}{Q_l^k-Q_j^k}=\frac{e^{i\beta_l}+p_l^k}{e^{i\beta_l}-e^{i\beta_j}+p_l^k-p_j^k}\\
=&\frac{e^{i\beta_l}+p_l^k}{(e^{i\beta_l}-e^{i\beta_j})(1+(p_l^k-p_j^k)/(e^{i\beta_l}-e^{i\beta_j}))}\\
=&\frac{e^{i\beta_l}}{e^{i\beta_l}-e^{i\beta_j}}+\frac{p_l^k}{e^{i\beta_l}-e^{i\beta_j}}-
\frac{e^{i\beta_l}}{(e^{i\beta_l}-e^{i\beta_j})^2}(p_l^k-p_j^k)+O(\sigma_k^2)\\
=&\frac{e^{i\beta_l}}{e^{i\beta_l}-e^{i\beta_j}}+\frac{e^{i\beta_l}p_j^k-e^{i\beta_j}p_l^k}{(e^{i\beta_l}-e^{i\beta_j})^2}+O(\sigma_k^2).
\end{align*}

Using
\begin{equation}\label{e-N}
N=2\sum_{j=0,j\neq l}^N\frac{e^{i\beta_l}}{e^{i\beta_l}-e^{i\beta_j}},
\end{equation}
we write (\ref{pi-N}) as
\begin{equation}\label{var-1}
\sum_{j=0,j\neq l}^N\frac{e^{i\beta_l}p_j^k-e^{i\beta_j}p_l^k}{(e^{i\beta_l}-e^{i\beta_j})^2}=E_1+O(\sigma_k^2)
\end{equation}
for $l=0,1,2,....,N$.
For convenience we set
$$p_l^k=e^{i\beta_l}m_l^k \quad \mbox{ and  }\quad \beta_{jl}=\beta_j-\beta_l $$
to  reduce (\ref{var-1}) to
\begin{align}\label{main-1}
&\sum_{j=0,j\neq l}^N\frac{e^{i\beta_{jl}}m_j^k}{(1-e^{i\beta_{jl}})^2}-\bigg (\sum_{j=0,j\neq l}^N\frac{e^{i\beta_{jl}}}{(1-e^{i\beta_{jl}})^2}\bigg )m_l^k\\
&=E+O(\sigma_k^2)+O(\delta_k \sigma_k) \nonumber
\end{align}
for $l=0,1.....,N$.
It is easy to verify that
\begin{equation}\label{trig-i-1}
\frac{e^{i\theta}}{(1-e^{i\theta})^2}=\frac{1}{2(\cos \theta -1)}=(-\frac 14)\frac{1}{\sin^2(\theta/2)}.
\end{equation}
To deal with coefficients of $m_j^k$ in (\ref{main-1}) we set
$$d_j=\frac{1}{\sin^2(\frac{j\pi}{N+1})},\quad j=1,...,N $$
and
$$D=\sum_{j=0,j\neq l}^N d_{|j-l|}. $$
Since $d_l=d_{N+1-l}$ it is easy to check that $D$ does not depend on $l$:
\begin{equation}\label{D-guess}
D=\sum_{k=1}^N d_k=\sum_{k=1}^N\frac{1}{\sin^2(\frac{k\pi}{N+1})}=\frac{N^2+2N}{3}.
\end{equation}
Now (\ref{main-1}) can be written as
\begin{equation}\label{main-2}
-\sum_{j\neq l,j=0}^N d_{|j-l|}m_j^k+Dm_l^k=E_1+O(\sigma_k^2),\quad l=0,....,N.
\end{equation}

For $l=0$, we have $\beta_0=0$ and $m_0^k=0$. Thus from (\ref{main-2}) we have
\begin{equation}\label{main-3}
-\sum_{j=1}^N d_jm_j^k =E_1+O(\sigma_k^2).
\end{equation}
If we take $(m^k_1,...,m^k_n)$ as unknowns in (\ref{main-2}), the last $N$ equations of (\ref{main-2}) ( for $l=1,...,N$) can be written as
\begin{equation}\label{e-m}
A\left(\begin{array}{c}
m^k_1\\
m^k_2\\
\vdots\\
m^k_N
\end{array}
\right)=E_1+O(\sigma_k^2).
\end{equation}
where
$$A=\left(\begin{array}{cccc}
D & -d_1 & ... & -d_{N-1} \\
-d_1 & D & ... & -d_{N-2} \\
\vdots & \vdots & ... & \vdots \\
-d_{N-1} & -d_{N-2} & ... & D
\end{array}
\right ) $$

Since $D=|d_1|+...+|d_N|$ and each $d_i>0$, we see that the matrix is invertible, thus $|m^k_i|=E$ for all $i$.
 
One consequence about $v_{0,k}$ is that the local maximum points of $v_{0,k}$, being perturbed by the amount of $O(\delta_k^*)$, are $O(e^{-\mu_k}\delta_k^*)$ away from the corresponding locations of the local maximum points of $v_k$. The reason is $v_k$ is non-degenerate near each blowup point.

\end{document}